\documentstyle[amstex,amscd]{article}
\title{Deformation theory of singular symplectic n-folds}
\author{Yoshinori Namikawa}
\date{ }
\begin{document}
\maketitle
\begin{center}
{\bf Introduction}
\end{center}

By a symplectic manifold (or a symplectic n-fold) 
we mean a compact Kaehler manifold of even 
dimension $n$ 
with a non-degenerate holomorphic 2-form $\omega$, 
i.e. $\omega^{n/2} $ is a nowhere-vanishing $n$-form. 
This notion is generalized to a variety with 
singularities. 
We call $X$ a projective symplectic variety 
if $X$ is a normal projective variety with 
rational Gorenstein singularities and 
if the regular locus $U$ of $X$ admits 
a non-degenerate holomorphic 2-form $\omega$.
 A symplectic variety will play an important role
 together with a singular Calabi-Yau variety 
in the generalized Bogomolov decomposition conjecture. 
Now that essentially a few examples of symplectic 
manifolds are discovered, it seems an important 
task to seek new symplectic manifolds by 
deforming symplectic varieties. 
In this paper we shall study a projective symplectic 
variety from a view point of deformation theory. 
If $X$ has a resolution $\pi : \tilde X \to X$ 
such that $(\tilde X, \pi^*\omega)$ is a 
symplectic manifold, we say that $X$ has a 
symplectic resolution. 
Our first results are concerned with a 
birational contraction map of a symplectic 
manifold. 
\vspace{0.2cm}

{\bf Proposition (1.4)}. {\em Let 
$\pi : \tilde X \rightarrow X$ be a birational 
projective morphism from a projective 
symplectic n-fold 
$\tilde X$ to a normal n-fold $X$.
 Let $S_i$ be the set of points 
$p \in X$ such that $\dim \pi^{-1}(p) = i$. 
Then $\dim S_i \leq n - 2i$. In particular, 
$\dim \pi^{-1}(p) \leq n/2$}.  \vspace{0.12cm}

{\bf Proposition (1.6)}. {\em  Let 
$\pi : \tilde X \rightarrow X$ be a birational 
projective morphism from a 
projective symplectic n-fold 
$\tilde X$ to a normal n-fold $X$. Then $X$ has 
only canonical singularities and its dissident 
locus $\Sigma_0$ has codimension at least 4 in $X$. 
Moreover, if $\Sigma \setminus \Sigma_0$ is non-empty, 
then $\Sigma\setminus\Sigma_0$ is a disjoint union 
of smooth varieties of dim $n-2$ with everywhere 
non-degenerate 2-forms. }  \vspace{0.12cm}

When $X$ has only an isolated singularity $p \in X$, 
every irreducible component of $\pi^{-1}(p)$ 
is Lagrangian (Proposition (1.11)).  
 In this situation it is conjectured that 
the exceptional locus is isomorphic 
to ${\bold P}^{n/2}$ 
with normal bundle $\Omega^1_{{\bold P}^{n/2}}$. 
Similar results to (1.4) and (1.11) are obtained 
independently by Wierzba [Wi].  

We shall exhibit four examples of birational 
contraction maps of symplectic 4-folds in (1.7). 
The second example shows that the Kaehler 
(projective) assumption of 
a symplectic manifold is not necessarily preserved 
under an elementary transformation. 
The fourth example deals with a symplectic manifold 
obtained as a resolution of certain quotient of 
a Fano scheme of lines on a cubic 4-fold.
As for a fiber space structure of a 
symplectic n-fold, see [Ma].  

After we study the birational contraction map 
of a symplectic manifold in section 1, we shall 
prove our main theorem in section 2: \vspace{0.15cm}

{\bf Theorem (2.2)}. {\em Let  $\pi : \tilde X \to X$ 
be a symplectic resolution of a projective 
symplectic variety $X$ of dimension $n$. 
Then the Kuranishi spaces ${\mathrm{Def}}(\tilde X)$ 
and ${\mathrm{Def}}(X)$ are both smooth 
of the same dimension. 
There exists a natural map $\pi_* : 
{\mathrm{Def}}(\tilde X) 
\to {\mathrm{Def}}(X)$ and $\pi_*$ is a finite covering 
\footnote{Precisely, there are open subsets 
$0 \in V \subset {\mathrm{Def}}(\tilde X)$ and 
$0 \in W \subset {\mathrm{Def}}(X)$ such that 
$\pi_*\vert_V : V \to W$ is a proper surjective 
map with finite fibers.}. 
Moreover, $X$ has a flat deformation to a smooth 
symplectic n-fold $X_t$. Any smoothing $X_t$ of $X$ 
is a symplectic n-fold obtained as a 
flat deformation of $\tilde X$.}
 \vspace{0.2cm}

(2.2) was proved by Burns-Wahl [B-W] 
for K3 surfaces. Given a one-parameter flat deformation 
$f: \cal X \to \Delta$ of such $X$ as (2.2), 
by Theorem, we could have a simultaneous 
resolution $\nu : \tilde{\cal X} 
\to {\cal X}'$ after a suitable finite base change 
${\cal X}' \to \Delta'$ of $\cal X$ by $\Delta' 
\to \Delta$.         

The same situation as (2.2) naturally arises for 
Calabi-Yau 3-folds; but the results for them 
are very partial as compared with symplectic case 
(cf. Example (2.4)).

On the other hand, it is natural to consider a 
symplectic variety which does not have a symplectic 
resolution; for example, such varieties appear in a 
work of O'Grady [O] as the moduli spaces of rank 2 
semi-stable sheaves on a K3 surface with $c_1 = 0$ 
and with even $c_2 \geq 6$. 
At the moment it is not clear when these varieties 
have flat deformations to symplectic manifolds. 
But we can prove that such varieties have 
unobstructed deformations: \vspace{0.12cm}

{\bf Theorem (2.5)}. {\em Let $X$ be a projective 
symplectic variety. 
Let $\Sigma \subset X$ be the singular locus.  
Assume that $\mathrm{codim}(\Sigma \subset X) \geq 4$. 
Then ${\mathrm{Def}}(X)$ is smooth. }  \vspace{0.15cm}

We shall give a rough sketch of the proof 
of Theorem (2.2) in the remainder. 

First note that $X$ has only rational 
Gorenstein singularities. 
Then the existence of the map 
$\pi_*$ follows from the 
fact that $R^1\pi_*{\cal O}_{\tilde X} = 0$ 
(cf. [Ko-Mo, (11.4)]). 

Let $U$ be the complement of $\Sigma_0$ in $X$
 and write $\tilde U$ for $\pi^{-1}(U)$. 
By (1.4) and (1.6), we can prove, 
roughly speaking, that a deformation of $\tilde X$ 
(resp. $X$) is equivalent to that of $\tilde U$ 
(resp. $U$). (See Proposition (2.1).) 
From this fact it follows that 
$\pi_* : {\mathrm{Def}}(\tilde X) \to 
{\mathrm{Def}}(X)$ 
is finite. 

Finally we compare the dimensions of tangent spaces of 
${\mathrm{Def}}(X)$ and ${\mathrm{Def}}(\tilde X)$ 
at the origin and then conclude 
that ${\mathrm{Def}}(X)$ 
is smooth. Since ${\mathrm{Def}}(\tilde X)$ 
is smooth by 
Bogomolov [Bo], we only have to prove that 
$\dim {\bold T}^1_X = \dim {\bold T}^1_U$ 
is not larger than 
$h^1(\tilde X, \Theta_{\tilde X}) = 
h^1(\tilde U, \Theta_{\tilde U})$.  
We need here a detailed description of the sheaf 
$T^1_U := 
{\underline{Ext}}^1({\Omega^1}_U, {\cal O}_U)$ 
(Lemma (1.9), Corollary (1.10)). 

The last statement will be proved in the following way. 
By the existence of a non-degenerate 2-form $\omega$, 
there is an obstruction to extending a holomorphic 
curve on $\tilde X$ sideways in a given one-parameter 
small deformation $\tilde{\cal X} \to \Delta^1$. 
Therefore, if we take a general curve of 
${\mathrm{Def}}(\tilde X)$ passing through the origin 
and take a corresponding small deformation of 
$\tilde X$, then no holomorphic curves survive 
(cf. [Fu, Theorem (4.8)]).  

Let $t \in {\mathrm{Def}}(X)$ be a generic point 
(that is, $t$ is outside the union of a countable number 
of proper subvarieties of ${\mathrm{Def}}(X)$).  
Since 
$\pi_* : {\mathrm{Def}}(\tilde X) \to {\mathrm{Def}}(X)$ 
is a finite covering, we may assume that 
$X_t$ has a symplectic resolution 
$\pi_t : {\tilde X}_t \to X_t$. By the argument above, 
${\tilde X}_t$ contains no curves. 
By Chow lemma [Hi], there is a bimeromorphic 
projective map $h : W \to X_t$ such that $h$ is factored 
through $\pi_t$. Since $h^{-1}(p)$ is 
the union of projective 
varieties for any point $p \in X_t$, ${\pi_t}^{-1}(p)$ 
is the union of Moishezon varieties. 
If $\pi_t$ is not an isomorphism, ${\pi_t}^{-1}(p)$ 
has positive dimension for some $p \in X_t$; 
hence ${\tilde X}_t$ contains curves, 
which is a contradiction. Thus $\pi_t$ 
is an isomorphism and $X_t$ is a (smooth) 
symplectic $n$-fold.         
\vspace{0.2cm}

The author thanks A. Fujiki for giving him 
invaluable informations on this topic. The first 
version of this paper was written in 1998. After 
that the author was informed that Wierzba [Wi] 
independently obtained similar results to 
(1.4) and (1.11).  \vspace{0.3cm}
  
\begin{center}
{\bf 1. Birational contraction maps of symplectic n-folds}
\end{center}

A symplectic $n$-fold means a symplectic {\em manifold} 
of dimension $n$.
We shall state three lemmas which will be used later. 
The first lemma is essentially a linear algebra. 
\vspace{0.15cm}

{\bf Lemma (1.1)}. {\em Let $V$ be a complex manifold with 
$\dim V = 2r$ and let $\omega$ be an everywhere 
non-degenerate holomorphic 2-form on $V$ (i.e. 
$\wedge^r \omega$ is nowhere-vanishing.)  Let $E$ 
be a subvariety of $V$ with $\dim E > r$. Then 
$\omega\vert_E$ is a non-zero 2-form on $E$. }
\vspace{0.15cm}

{\bf Lemma (1.2)}. {\em Let $f: V \to W$ be a birational 
projective morphism from a complex manifold $V$ 
to a normal variety $W$. Let $p \in W$ and assume that 
the germ $(W, p)$ of $W$ at 
$p$ has rational singularities. 
Assume that $E := f^{-1}(p)$ is a simple normal crossing 
divisor of $V$. Then $H^0(E, \hat\Omega^i_E) = 0$ for 
$i > 0$, where $\hat\Omega^i_E 
:= \Omega^i_E /(torsion)$.} \vspace{0.12cm}

{\em Proof}. Denote by $F^{\cdot}$ (resp. $W_{\cdot}$) 
the Hodge filtration (resp. weight filtration) of 
$H^i(E): = H^i(E, \bold C)$. 
Note that these two filtrations 
give a mixed Hodge structure on $H^i(E)$. 
Since $E$ is a proper algebraic scheme, 
$Gr^W_j(H^i(E)) = 0$ for $j > i$. 

 Assume that $H^0(E, \hat{\Omega}^i_E) 
\neq 0$. Then $Gr^i_FGr^W_i(H^i(E)) \neq 0$. 
By the Hodge symmetry 
$Gr^0_FGr^W_i(H^i(E)) \neq 0$, 
and hence $Gr^0_F(H^i(E)) = 
H^i(E, {\cal O}_E) \neq 0$.  

On the other hand, $(R^if_*{\cal O}_V)_p = 0$ 
for $i > 0$  
because $(W, p)$ has only rational singularities. 
Take a sufficiently small open neighborhood $V'$ of 
$f^{-1}(p)$ in $V$. There exists a commutative diagram 
of Hodge spectral sequences

\begin{equation}
\begin{CD}
H^k(V', \Omega^j_{V'})  @=> H^{j+k}(V', {\bold C}) \\
@VVV @VVV \\
H^k(E, \hat\Omega^j_E) @=>  
H^{j+k}(E, {\bold C})  
\end{CD}
\end{equation}

Note that 
$H^i(V', {\bold C}) \cong H^2(E, {\bold C})$. 
Denote by $F_1^{\cdot}$ (resp. $F_2^{\cdot}$) 
the filtrations on $H^i(V', {\bold C})$  
(resp. $H^i(E, {\bold C})$) 
induced by the spectral sequences. 
There is a natural surjection 
$Gr^0_{F_1}H^i(V', {\bold C}) \to 
Gr^0_{F_2}H^i(E, {\bold C})$.  As 
$(R^if_*{\cal O}_V)_0 = 0$ for $i > 0$, 
$H^i(V', {\cal O}_{V'}) = 0$. Therefore we have 
$Gr^0_{F_1}H^2(V', {\bold C}) = 
Gr^0_{F_2}H^i(E, {\bold C}) = 0$. 
Since the second spectral sequence degenerates 
at $E_1$ terms, 
$H^i(E, {\cal O}_E) = 0$, 
which is a contradiction. \vspace{0.12cm}

{\bf Lemma(1.3)}. {\em Let $V$ be a 
symplectic n-fold and let $H$ be a 
smooth 3-dimensional 
subvariety of $V$ containing a smooth 
rational curve $C$ with 
$N_{H/V}\vert_C \cong {\cal O}^{\oplus n-3}$. 
Assume that 
$N_{C/H} \cong {\cal O}(-1)\oplus{\cal O}(-1)$, 
${\cal O}(-2)\oplus{\cal O}$ or 
${\cal O}(-3)\oplus{\cal O}(1)$. 
Then ${\mathrm{Hilb}}(V)$ is smooth of dimension 
(n-2) at [C].
Moreover, in this case, 
$N_{C/V} \cong 
{\cal O}^{\oplus (n-2)}\oplus{\cal O}(-2)$ or  
${\cal O}^{\oplus (n-4)}\oplus{\cal O}(-1)
\oplus{\cal O}(1)\oplus{\cal O}(-2)$}. \vspace{0.15cm}

{\em Proof}. We shall prove that the Hilbert 
scheme (functor) has the $T^1$-lifting property at [C]; 
then ${\mathrm{Hilb}}(V)$ 
is smooth at [C] by [Ra, Ka 1]. 
Let $S_m$ be the spectrum of the Artinian ring 
$A_m = {\bold C}[t]/(t^{m+1})$. 
Set ${V}_m := V \times_{S_0}S_m$. 
Let $C_m \subset {V}_m$ be an infinitesimal 
displacement of $C$ to m-th order, and let 
$C_{m-1} := 
C_m \times_{{V}_m}{V}_{m-1}$. 
We have to prove that 

$$  H^0(N_{C_m/{V}_m}) \rightarrow 
H^0(N_{C_{m-1}/{V}_{m-1}}) $$

is surjective. Let $\omega$ be a non-degenerate 
2-form on $V$. Then $\omega$ lifts to 
an element 
$\omega_m \in H^0(\Omega^2_{{V}_m/S_m})$ 
in such a way that 
$\wedge^{n/2} \omega_m 
\in H^0(\Omega^n_{{V}_m/S_m})$ 
is a nowhere vanishing section (that is, 
if we identify 
$H^0(\Omega^n_{{V}_m/S_m})$ with $A_m$, 
then $\wedge^{n/2} \omega_m$ corresponds 
to an invertible emement of $A_m$). 
The 2-form $\omega_m$ induces a pairing 

$$  \Theta_{{V}_m/S_m}\vert_{C_m} 
\times \Theta_{{V}_m/S_m}\vert_{C_m} 
\rightarrow {\cal O}_{C_m}. $$

Since this pairing vanishes on 
$\Theta_{C_m/S_m} \times \Theta_{C_m/S_m}$ 
and since $\omega_m$ is non-degenerate, 
one has a surjection 

$$  \alpha_m : N_{C_m/{V}_m} 
\rightarrow \Omega^1_{C_m/S_m} $$

by the exact sequence 

$$  0 \to \Theta_{C_m/S_m} \to 
\Theta_{{V}_m/S_m}\vert_{C_m} 
\to N_{C_m/{V}_m} \to 0. $$

Let us first consider the case when $m = 0$. 
By assumption, we have 
$N_{H/V}\vert_C \cong {\cal O}^{\oplus n-3}$; 
hence by the exact sequence 

$$   0  \to N_{C/H} \to N_{C/V} 
\to N_{H/V}\vert_C \to 0 $$

we see that $N_{C/V}$ is isomorphic to 
${\cal O}(-1)
\oplus{\cal O}(-1)\oplus{\cal O}^{\oplus (n-3)}$, 
${\cal O}(-3)\oplus{\cal O}(1)
\oplus{\cal O}^{\oplus (n-3)}$, 
${\cal O}(-2)\oplus{\cal O}(-1)
\oplus{\cal O}(1)\oplus{\cal O}^{\oplus (n-4)}$ 
or ${\cal O}(-2)\oplus{\cal O}^{\oplus (n-2)}$. 
By the existence of the surjection $\alpha_0$, 
the first two cases are excluded. In particular, 
it is checked that $H^1(\mathrm{Ker}(\alpha_0)) = 0$. 
Note that this implies that 
$H^1(\mathrm{Ker}(\alpha_m)) = 0$ for all 
$m$ because there are exact sequences 
$0 \to \mathrm{Ker}(\alpha_0) \to 
\mathrm{Ker}(\alpha_m) \to \mathrm{Ker}(\alpha_{m-1}) 
\to 0$.   

Next consider the following commutative diagram 
with exact columns and exact rows 

\begin{equation}
\begin{CD}
@. 0 @. 0 @. 0  \\
@. @VVV @VVV @VVV \\
0 @>>> H^0(\mathrm{Ker}(\alpha_0)) @>>> 
H^0(\mathrm{Ker}(\alpha_m)) @>{\psi_m}>> 
H^0(\mathrm{Ker}(\alpha_{m-1})) \\
@. @VVV @VVV @VVV \\
0 @>>> H^0(N_{C/{V}}) @>>> 
H^0(N_{C_m/{V}_m}) @>{\phi_m}>> 
H^0(N_{C_{m-1}/{V}_{m-1}}) \\
@. @V{\tau_0}VV @V{\tau_m}VV @V{\tau_{m-1}}VV \\
0 @>>> H^0(\Omega^1_C) @>>> H^0(\Omega^1_{C_m/S_m}) 
@>{\varphi_m}>> H^0(\Omega^1_{C_{m-1}/S_{m-1}})   
\end{CD}
\end{equation}

As we remarked above, 
$H^1(\mathrm{Ker}(\alpha_m)) = 0$ for all $m$, 
hence $\tau_0$, $\tau_m$ and $\tau_{m-1}$ are 
surjective. By the same reason, $\psi_m$ 
is also surjective. By the Hodge theory, 
$\varphi_m$ is surjective. 
It now follows from the diagram above that 
$\phi_m$ is surjective. Thus the Hilbert scheme 
${\mathrm{Hilb}}(V)$ is smooth at [C]. 
Its dimension equals $h^0(N_{C/V}) = n-2$.    

{\bf Proposition (1.4)}. {\em Let 
$\pi : \tilde X \rightarrow X$ be a birational 
projective morphism from a projective symplectic n-fold 
$\tilde X$ to a normal n-fold $X$.
 Let $S_i$ be the set of points 
$p \in X$ such that $\dim \pi^{-1}(p) = i$. 
Then $\dim S_i \leq n - 2i$. In particular, 
$\dim \pi^{-1}(p) \leq n/2$}.  \vspace{0.15cm}

{\em Proof}.  For a non-empty $S_i$, we take 
an irreducible component $R_i$ of 
$\pi^{-1}(S_i)$ in such a way that 

(1) by $\pi$, $R_i$ dominates an irreducible component 
of $S_i$ with $\dim S_i$, and 

(2) a general fiber of $R_i \to \pi (R_i)$ has dimension 
$i$. 
  
Put $l := n - \dim R_i$. By definition 
$\dim S_i = n - l - i$. We shall prove that 
$\dim S_i \geq n - 2l$; if this holds, then 
$i \leq l$, and hence $\dim S_i \leq n - 2i$. 
When $l \geq n/2$, then clearly 
$\dim S_i \geq n - 2l$. 
We assume that $l < n/2$.
We shall derive a contradiction 
assuming that $\dim S_i < n - 2l$ and 
assuming that $S_i$ is irreducible.  
When $S_i$ is not irreducible, 
it is enough only to replace $S_i$ by $\pi(R_i)$.  

Take a birational projective morphism 
$\nu : Y \to \tilde X$ in such a way that 
$F := (\pi\circ\nu)^{-1}(S_i)$ becomes a 
divisor of a smooth n-fold $Y$ with 
normal crossings. Set $f = \pi\circ\nu$.

A non-degenerate 2-form $\omega$ on $\tilde X$ 
is restricted to a non-zero 2-form on 
$R_i$ because $\dim R_i > n/2$ (Lemma (1.1)).   
Therefore we have a non-zero element 
$\nu^*\omega\vert_{F} \in 
H^0(F, \hat{\Omega}^2_{F})$.

For a general point $p \in S_i$, the fiber 
$F_p$ of the map $F \to S_i$ is a normal crossing 
variety. Hence, if we take a suitable open set 
$U_i$ of $S_i$ and replace $F$ by 
$(\pi\circ\nu)^{-1}(U_i)$, then
the sheaf 
${\hat\Omega}^2_{F}$ has a filtration 
$f^*{\Omega}^2_{S_i} \subset 
{\cal F} \subset {\hat\Omega}^2_{F}$ 
with the exact sequences 

$$  0 \to {\cal F} \to {\hat\Omega}^2_{F} 
\to {\hat\Omega}^2_{F/S_i} \to 0 $$ 

$$  0 \to f^*{\Omega}^2_{S_i} \to {\cal F} 
\to f^*{\Omega}^1_{S_i}\otimes{\hat\Omega}^1_{F/S_i} 
\to 0 $$

Let us prove that the 2-form 
$\nu^*\omega\vert_{F}$ is not the 
pull-back of any 2-form on $S_i$. 
 Write $n = 2r$. 
Assume that $\omega' := \nu^*\omega\vert_{F}$ 
is the pull-back 
of a 2-form on $S_i$. Then 
$\wedge^{r-l}\omega' = 0$ because 
$\dim S_i < 2r - 2l$.  
On the other hand, take a general point $q \in R_i$. 
Since $R_i$ is a submanifold of $\tilde X$ 
of codimension $l$ around $q$, it is checked by 
linear algebra that 
$\wedge^{r-l}(\omega\vert_{R_i}) \ne 0$ 
in a open neighborhood of $q \in R_i$. 
Let $F'$ be the union of irreducible 
components of $F$ 
which dominate $R_i$ by $\nu$. Since 
$\nu^*\omega\vert_{F'} = 
(\nu\vert_{F'})^*(\omega\vert_{R_i})$, 
$\wedge^{r-l}(\nu^*\omega\vert_{F'}) \ne 0$. 
Since $(\wedge^{r-l}\omega')\vert_{F'} = 
\wedge^{r-l}(\nu^*\omega\vert_{F'})$, 
this implies that $\wedge^{r-l}\omega' \ne 0$, 
which is a contradiction.      \vspace{0.15cm}

Let $F_p$ be the fiber of $F \to S_i$ over 
$p \in S_i$. Note that $F_p$ is a 
normal crossing variety for a general point 
$p \in S_i$. Then, by the exact sequence, 
one can see that 
$H^0(F_p, \hat{\Omega}^1_{F_p}) \neq 0$ or 
$H^0(F_p, \hat{\Omega}^2_{F_p}) \neq 0$. 

Take an $l+i$ dimensional complete intersection 
$H = H_1 \cap H_2 ... \cap H_{n-l-i}$ 
of very ample divisors of $X$ passing 
through a general point $p \in S_i$. 
(When $l + i = n$, we put $H = X$.) 
Then $H$ has only rational singularities. 
Put $\tilde H := f^{-1}(H)$ and put 
$g := f\vert_{\tilde H}$.
Note that $g^{-1}(p) = F_p$ is a divisor of 
$\tilde H$ with normal crossings. 
By Lemma (1.2) 
$H^0(F_p, \hat{\Omega}^i_{F_p}) = 0$ 
for $i > 0$, which is a contradiction.
\vspace{0.15cm}

{\bf Corollary (1.5)}. {\em  Let 
$\pi : \tilde X \rightarrow X$ be a birational 
projective morphism from a projective symplectic n-fold 
$\tilde X$ to a normal n-fold $X$. 
Then any $\pi$-exceptional divisor is mapped 
onto an $(n-2)$-dimensional subvariety of $X$ 
by $\pi$. } \vspace{0.15cm}

{\em Proof}. Take a $\pi$-exceptional divisor $E$.  
For some $i \geq 1$ we can take the $E$ as an 
$R_i$ in the proof of Proposition (1.4). Then 
$\dim \pi(E) \geq n - 2$. \vspace{0.2cm}
 
Let $X$ be a normal variety of dim $n$ with canonical 
singularities. Let $\Sigma$ be the singular 
locus of $X$. By [Re] there is a closed subset 
$\Sigma_0 \subset \Sigma$ such that each point 
of $\Sigma \setminus \Sigma_0$ has an analytic 
open neighborhood in $X$ isomorphic to 
(rational double point) $\times ({\bold C}^{n-2}, 0)$. 
The locus $\Sigma_0$ is called the dissident locus. 
Generally we have $\dim \Sigma_0 \leq n-3$. 
But, when $X$ has a symplectic resolution, we have 
a stronger result.  
\vspace{0.15cm} 

{\bf Proposition (1.6)}. {\em  Let 
$\pi : \tilde X \rightarrow X$ be a birational 
projective morphism from a projective symplectic n-fold 
$\tilde X$ to a normal n-fold $X$. Then $X$ has 
only canonical singularities and its dissident 
locus $\Sigma_0$ has codimension at least 4 in $X$. 
Moreover, if $\Sigma \setminus \Sigma_0$ is non-empty, 
then $\Sigma\setminus\Sigma_0$ is a disjoint union 
of smooth varieties of dim $n-2$ with everywhere 
non-degenerate 2-forms. } 
\vspace{0.15cm}

{\em Proof}: 

{\bf (1.6.1)} {\em $\Sigma$ has no (n-3)-dimensional 
irreducible components.} \vspace{0.15cm}

We shall derive a 
contradiction by assuming that $\Sigma$ has 
an $(n-3)$-dimensional irreducible component. 
Let $H := H_1 \cap H_2 \cap ... \cap H_{n-3}$ 
be a complete intersection of very ample 
divisors of $X$. The $H$ intersects 
the $(n-3)$-dimensional component in finite points. 
Let $p \in H$ be one of such points. 
Let $H' := \pi^{-1}(H)$. 
Since there are no exceptional divisors of 
$\pi$ lying on the $(n-3)$-dimensional component 
of $\Sigma$, $\pi\vert_{H} : H' \to H$ 
gives a small resolution of $H$ around $p$. 
Pick an irreducible curve $C$ from 
$\pi\vert_{H'}^{-1}(p)$. The $C$ is isomorphic to 
${\bf P}^1$, and its normal bundle 
$N_{C/H'}$ in $H'$ is isomorphic to one of three vector 
bundles ${\cal O}(-1)\oplus{\cal O}(-1)$, 
${\cal O}(-2)\oplus{\cal O}$ or 
${\cal O}(-3)\oplus{\cal O}(1)$. 
Note that the Hilbert scheme ${\mathrm{Hilb}}(\tilde X)$ 
has at most dimension $(n-3)$ because $C$ can only move 
in $\tilde X$ along the $(n-3)$-dimensional component 
of $\Sigma$. However this contradicts Lemma (1.3).
\vspace{0.15cm}

{\bf (1.6.2)}. By (1.6.1) we only 
have to observe the irreducible components 
of $\Sigma$ with dimension $n-2$. So we 
replace $\Sigma$ by an irreducible component 
of $\Sigma$ with dim $n-2$.
We shall derive a contradiction by assuming 
that $\dim \Sigma_0 = n-3$. 

Let $H := H_1 \cap H_2 \cap ... \cap H_{n-3}$ 
be a complete intersection of very ample divisors of 
$X$. Then $\tilde H := \pi^{-1}(H)$ is a crepant 
resolution of $H$. Set $\Lambda := \Sigma \cap H$ 
and $\Lambda_0 := \Sigma_0 \cap H$. 
Note that $\Lambda_0$ consists of finite points. 
Write $\tau : \tilde H \to H$ for the restriction 
$\pi\vert_{\tilde H}$ of $\pi$ to $\tilde H$. 
Every fiber of $\tau$ has at most 
dimension one because, if some fibers are 
2-dimensional, then there is a prime divisor of 
$\tilde X$ lying on $\Sigma_0$, 
which contradicts Corollary (1.5). 

We shall show that $\Lambda$ is a smooth curve 
and that ${\mathrm{Exc}}(\tau)$ is locally isomorphic 
to the product of $\Lambda$ and a tree of 
${\bold P}^1$'s. If so, then $H$ must have 
rational double points of the same type along 
$\Lambda$ and this is a contradiction. 
A contradiction  will be deduced in several steps.

{\bf (i)} Take a point $p_0 \in \Lambda_0$. 
We only have to argue locally around $p_0$. 
Since $H$ has rational singularities and 
since $\tau^{-1}(p_0)$ is 1-dimensional, 
$\tau^{-1}(p_0)$ is a tree of ${\bold P}^1$'s. 
Let $C_1, ..., C_m$ be the irreducible 
components of $\tau^{-1}(p_0)$. 
Let us compute the normal bundle 
$N_{C_i/\tilde H}$. 
Take a sufficiently small open neighborhood $U$ of 
$\cup C_i \subset \tilde H$. 
Since $H$ has only rational singularities, 
we have $H^1(U, {\cal O}_U) = 0$. 
Let $I_i$ be the defining ideal of $C_i$ in $U$. 
Then, by the exact sequence 
$H^1(U, {\cal O}_U) \to H^1(C_i, {\cal O}_U/I^2_i) 
\to H^2(U, I^2_i) = 0$ we know that 
$H^1(C_i, {\cal O}_U/I^2_i) = 0$. 
By another exact sequence 
$H^0(C_i, {\cal O}_U/I^2_i) \to 
H^0(C_i, {\cal O}_U/I_i)(= {\bold C}) \to 
H^1(C_i, I_i/I^2_i) \to 
H^1(C_i, {\cal O}_U/I^2_i) = 0$, 
we know that $H^1(C_i, I_i/I^2_i) = 0$ 
because the first map is surjective. 
Therefore, $N_{C_i/\tilde H} \cong 
{\cal O}(-1)\oplus{\cal O}(-1)$,
 ${\cal O}(-2)\oplus{\cal O}$ 
or 
${\cal O}(-3)\oplus{\cal O}(1)$.

By Lemma (1.3), the Hilbert scheme 
${\mathrm{Hilb}}(\tilde X)$ is smooth of 
dimension $(n-2)$ at [$C_i$]. 
This fact tells us two things. 
\vspace{0.15cm}

{\bf (i-a)}: {\em Each $C_i$ moves inside $\tilde H$}; 
in fact, if $C_i$ is rigid in $\tilde H$, 
then ${\mathrm{Hilb}}(\tilde X)$ possibly has only (n-3) 
parameter at [$C_i$] corresponding to a 
displacement of $C_i \subset \tilde X$ along 
$\Sigma_0$, which is a contradiction. \vspace{0.15cm}
  
{\bf (i-b)}: {\em We have $N_{C_i/\tilde H} 
\cong {\cal O}\oplus{\cal O}(-2)$, in particular, 
$N_{C_i/\tilde X} \cong 
{\cal O}^{\oplus (n-2)}\oplus{\cal O}(-2)$.}
  
This fact can be proved by using Grothendieck's 
Hilbert scheme (cf. [Ko 1, Chap. I]):   
Let ${\mathrm{Hilb}}(\tilde X/X)$ be the relative Hilbert 
scheme for $\pi : \tilde X \to X$. 
Since $C_i$ is contained in a fiber of $\pi$, 
${\mathrm{Hilb}}(\tilde X)$ coincides 
with ${\mathrm{Hilb}}(\tilde X/X)$ 
at $[C_i]$. Therefore ${\mathrm{Hilb}}(\tilde X/X)$ 
is smooth of dimension (n-2) at $[C_i]$. 
Moreover, the irreducible component of 
${\mathrm{Hilb}}(\tilde X/X)$ containing 
$[C_i]$ dominates an 
$(n-2)$-dimensional irreducible component of $\Sigma$ 
by the map ${\mathrm{Hilb}}(\tilde X/X) \to X$. 
By the universal property of the 
relative Hilbert scheme, we have 
${\mathrm{Hilb}}(\tilde H/H) \cong 
{\mathrm{Hilb}}(\tilde X/X) \times_X H$, 
and hence ${\mathrm{Hilb}}(\tilde H/H)$ is smooth of 
dimension 1 at $[C_i]$ by Bertini theorem. 
Since ${\mathrm{Hilb}}(\tilde H)$ coincides with 
${\mathrm{Hilb}}(\tilde H/H)$ at $[C_i]$, 
this implies that ${\mathrm{Hilb}}(\tilde H)$ 
is smooth of dimension 1 at $[C_i]$. 
Therefore we have $N_{C_i/\tilde H} 
\cong {\cal O}\oplus{\cal O}(-2)$.  
\vspace{0.15cm}
    
{\bf (ii)} We shall prove that $\Lambda$ is irreducible 
around $p_0 \in \Lambda_0$.
By (i-a) there are no flopping curves in 
${\mathrm{Exc}}(\tau)$, 
hence $\tau$ is a {\em unique} crepant resolution 
of $H$. Therefore, we can construct $\tau$ 
locally around $p_0$ in the following manner. 
Let $\Lambda_1, ..., \Lambda_n$ be the irreducible 
components of $\Lambda$ at $p_0$. 
Blow up $H$ at first along the defining 
ideal $I_1$ of the reduced subscheme $\Lambda_1$ 
and take its normalization. 
We shall prove that $\tau$ is factorized by 
this composition of blow-up and normalization.
We shall argue along the line of 
[Re 2, \S 2.12-15].  First note that $H$ 
is a cDV point by [Re 1, Theorem (2.2)]. 
Let us view $H$ as a total space of a flat family 
of surface rational double points over a 
disc $\Delta^1$. 
The $\tau$ then can be viewed as a 
simultaneous (partial) resolution of 
this flat family. Let $F_1, ..., F_l$ 
be the irreducible components of 
${\mathrm{Exc}}(\tau)$ which dominate 
$\Lambda_1$. 
There is a unique positive divisor 
$F = \Sigma a_iF_i$ such that 
$F$ meets each general fiber 
$\tilde H_t$ ($t \in \Delta^1$) 
in the sum of the Artin's fundamental cycles 
for the rational double points 
$H_t \cap \Lambda_1$.
 Since there are no rigid $\tau$-exceptional 
curves, any $\tau$-exceptional curve 
$C$ moves along some $\Lambda_i$. 
If $C$ moves along $\Lambda_i$  with $i > 1$, 
then $(-F.C) = 0$. 
If $C$ moves along $\Lambda_1$, 
then $(-F.C) \geq 0$ by the definition of $F$.  
Therefore, $-F$ is $\tau$-nef divisor. 
At a general point of $\Lambda_1$, 
$\tau_*{\cal O}_{\tilde H}(-F)$ 
coincides with the defining ideal 
sheaf $I_1$ of the reduced subscheme $\Lambda_1$. 
Since every fiber of $\tau$ has dimension $\leq 1$, 
$\tau_*{\cal O}_{\tilde H}(-F) \cong I_1$ 
(cf. [Re 2, (2.14)].  
Since $-F$ is a $\tau$-nef, $\tau$-big divisor, 
the natural map 
$\tau^*\tau_*{\cal O}_{\tilde H}(-F) \to 
{\cal O}_{\tilde H}(-F)$ is surjective. 
Thus the ideal 
$\tau^{-1}I_1 \subset {\cal O}_{\tilde H}$ 
is invertible. Let $H'_1$ be the blowing up of 
$H$ along $I_1$. 
Then $\tau$ is factorized as 
$\tilde H \to H'_1 \to H$.       

Take an irreducible component of the singular 
locus of the resulting 3-fold 
which dominates $\Lambda_1$. 
Blow up the 3-fold along the 
defining ideal of this irreducible component 
with reduced structure, 
and then take the normalization. 
Repeating such procedure resolves singularities 
along general points of $\Lambda_1$. 
Denote by $\tau_1 : H_1 \to H$ the resulting 3-fold. 
Next take an irreducible component of 
${\mathrm{Sing}}(H_1)$ which dominates $\Lambda_2$. 
Blow up $H_1$ along the defining ideal of it and 
take the normalization. By repeating them, 
$\tau$ is finally decomposed as 

$$  \tilde H = H_n \stackrel{\tau_n}\to H_{n-1} 
\stackrel{\tau_{n-1}}\to ... 
\stackrel{\tau_1}\to H $$

We shall derive a contradiction by assuming 
$n \geq 2$. By (i-b) there is a smooth surface 
$E_i \subset \tilde H$ which has a 
${\bold P}^1$-bundle structure containing $C_i$ 
as a fiber. These surfaces $E_i$ are mapped onto the 
{\em same} irreducible component of $\Lambda$
 by $\tau$; indeed, if 
$C_i \cap C_j \neq \emptyset$ and 
$\tau(E_i) \neq \tau(E_j)$, then 
$E_i \cap E_j =$ \{one point\}, 
which is a contradiction because both 
$E_i$ and $E_j$ are Cartier divisors of $\tilde H$. 
Moreover, $\tau(E_i) \neq \Lambda_n$. 
Indeed, suppose to the contrary. 
Then $C_1, ...,C_m$ are all contracted 
to a point by $\tau_n$. At the same time, 
all exceptional divisors of $\tau$ lying on 
$\Lambda_1, ..., \Lambda_{n-1}$ are contracted 
to curves. By the construction of $\tau_i$'s, 
this is a contradiction.  

On the other hand, the decomposition of 
$\tau$ explained above depends on 
the ordering of the irreducible components 
of $\Lambda$. Thus we have 
$\tau(E_i) \neq \Lambda_k$ for any 
$k \geq 1$, which is obviously a 
contradiction. \vspace{0.15cm}

{\bf (iii)}  Let $E_i \subset \tilde H$ be a 
smooth divisor mentioned above.  
It has a ${\bold P}^1$-bundle structure 
containing $C_i$ as a fiber. 
We shall prove \vspace{0.15cm}

{\bf (iii-a)}: {\em ${\mathrm{Exc}}(\tau)$ is a 
divisor with simple normal crossings}; 
\vspace{0.13cm}

{\bf (iii-b)}: ${\mathrm{Exc}}(\tau) = 
\bigcup_{1 \leq i \leq m}E_i$; \vspace{0.13cm}

{\bf (iii-c)}: {\em If $C_i \cap C_j = \emptyset$, 
then $E_i \cap E_j = \emptyset$.
If $C_i \cap C_j \neq \emptyset$, then 
$E_i \cap E_j$ is a section of at least one 
of the ${\bold P}^1$-bundles $E_i$ and $E_j$}; 
\vspace{0.13cm}

First we shall prove (iii-b). Since $\Lambda$ 
is irreducible by (ii), $\tilde H = H_1$ and 
$\tau = \tau_1$ in the notation of (ii). 
The $\tau_1$ is decomposed into blowing ups 
along irreducible reduced centers 
(followed by normalizations):

$$  \tilde H \to ... \stackrel{\sigma_3}\to H^{(2)} 
\stackrel{\sigma_2}\to H^{(1)} 
\stackrel{\sigma_1}\to H $$

By the construction, ${\mathrm{Exc}}(\sigma_k)$ 
has a fibration over an irreducible curve whose 
general fiber is isomorphic to ${\bold P}^1$ or a 
reducible line pair. When a general fiber of the 
fibration is irreducible, the special fiber must 
be irreducible. Indeed, if the special fiber 
contains more than one irreducible component, 
then the proper transform of some of them to 
$\tilde H$ becomes a rigid rational curve, 
which contradicts (i-a). In this case 
${\mathrm{Exc}}(\sigma_k)$ is irreducible.

When a general fiber of the fibration is reducible, 
the special fiber must have one or two irreducible 
components because, if it has more than two 
irreducible components, then the proper transform 
of some of them to $\tilde H$ becomes a 
rigid rational curve. 

If the special fiber has exactly two 
irreducible components, then 
${\mathrm{Exc}}(\sigma_k)$ has exactly two 
irreducible components. 

We shall prove that if the special fiber is 
irreducible, then ${\mathrm{Exc}}(\sigma_k)$ 
is also irreducible. 
Suppose to the contrary. 
Denote by $C$ the special fiber. 
Then ${\mathrm{Exc}}(\sigma_k)$ has exactly 
two irreducible components $F$ and $F'$. 
Each of them has a fibration over an 
irreducible curve, and the special fiber 
moves (as a 1-cycle on $H^{(k)}$) in both 
$F$ and $F'$.  
Let $\tilde F$ (resp. $\tilde F'$) be the 
proper transform of $F$ (resp. $F'$) by 
$\tilde{\sigma} : \tilde H \to H^{(k)}$. 
The $\tilde F$ (resp. $\tilde F'$) has a 
fibration over an irreducible curve containing 
the proper transform $\tilde C$ of $C$ 
in a special fiber. 
The special fiber has only one irreducible component 
$\tilde C$  because if it contains more, 
then $\tilde C$ is a rigid rational curve 
\footnote{ By Theorem (2.2) 
from [Re 1], we know that $H^{(k-1)}$ has only cDV 
sinularities. Put $p := {\sigma_k}(C)$. 
The germ $(H^{(k-1)}, p)$ is  then isomorphic to a 
hypersurface singularity $x^2 + f(y,z,w) = 0$ 
with $deg(f) \geq 3$. There is an involution 
$\iota$ of $(H^{(k-1)}, p)$ defined by 
$x \to -x$, $y \to y$, $z \to z$ and $w \to w$. 
Since $\tilde H$ is a unique crepant resolution 
of $H^{(k-1)}$, the $\iota$ lifts to an involution 
$\tilde{\iota}$ of 
$({\tilde H}, \tilde{\sigma}^{-1}(C))$. 
By $\tilde{\iota}$, $\tilde F$ and $\tilde F'$ 
are interchanged. Therefore, if the special fiber 
for $\tilde F$ is reducible, then the special 
fiber for $\tilde F'$ is also reducible. 
Any $\tilde{\sigma}$-exceptional divisor does 
not contain $\tilde C$. Since $\tilde C$ 
does not move in $\tilde F$ or $\tilde F'$, 
$\tilde C$ must be rigid in $\tilde H$.}      
and this contradicts (i-a).

Thus $\tilde C$ moves (as a 1-cycle) in both 
$\tilde F$ and $\tilde F'$. 
On the other hand, since $\tilde C$ 
coincides with one of $C_i$'s, $\tilde C$ 
should move as fibers in only one smooth 
${\bold P}^1$-bundle by (i-b). This is a 
contradiction.  

As a consequence, we know that 
${\mathrm{Exc}}(\tau)$ has exactly $m$ irreducible 
components. Since $E_i$'s are contained in 
${\mathrm{Exc}}(\tau)$, (iii-b) holds.

We shall next prove (iii-c) and (iii-a). 
The first statement of (iii-c) is clear. 
Assume that $C_i \cap C_j \neq \emptyset$. 
Denote by $p_i : E_i \to \Delta^1$ 
(resp. $p_j : E_j \to \Delta^1$) the 
${\bold P}^1$-bundle structure of $E_i$ 
(resp. $E_j$) whose central fiber over 
$0 \in \Delta^1$ is $C_i$ (resp. $C_j$). 
The intersection $E_i \cap E_j$ is 
multi-sections of $p_i$ and $p_j$ of degree 
$n_i$ and $n_j$ respectively. Suppose that 
$n_i > 1$ and $n_j > 1$. 

Let ${\cal C}$ be the set of all irreducible 
curves on $\tilde H$ which are fibers of $p_i$ 
or $p_j$. {\em For $l$, $l' \in {\cal C}$, 
we say $l$ and $l'$ are equivalent if there is 
a sequence of the elements of ${\cal C}$: 
$l_0 := l$, $l_1$, ..., $l_{k_0-1}$, $l_{k_0} := l'$ 
such that $l_k \cap l_{k+1} \neq \emptyset$ 
for any $k$}. This is an equivalence relation of 
${\cal C}$.  

Take a general fiber $l^*$ of $p_i$ and consider 
the set ${\cal C}(l^*)$ of all curves which are 
equivalent to $l^*$. Note that ${\cal C}(l^*)$ 
is a finite set  consisting of smooth rational curves. 
Pick up an element $l \in {\cal C}(l^*)$ 
which is a fiber of $p_i$. Then there are at least 
$n_i$ fibers of $p_j$ which intersect $l$. 
Similarly, for any element $m \in {\cal C}(l^*)$ 
which is a fiber of $p_j$, there are at least 
$n_j$ fibers of $p_i$ which intersect $m$.
This implies that ${\cal C}(l^*)$ is {\em not} a 
tree of ${\bold P}^1$'s.

On the other hand, ${\cal C}(l^*)$ is contained 
in a fiber of $\tau : \tilde H \to H$, 
which is a contradiction. Therefore, $n_i = 1$ or 
$n_j = 1$. One can assume that $n_i = 1$. 
In this case, $E_i \cap E_j$ is a section of $p_i$, 
and $E_j$ intersects $E_i$ with multiplicity one 
along $E_i \cap E_j$ because, if not, then it 
contradicts the fact that each fiber of $\tau$ is a 
tree of ${\bold P}^1$'s.  Since there are no triple 
points in $E_1 \cup E_2 \cup ... \cup E_m$, 
(iii-a) and (iii-c) hold.  \vspace{0.15cm}

{\bf (iv)} We shall prove that $E_i \cap E_j$ is sections 
of {\em both} ${\bold P}^1$-bundles $E_i$ and $E_j$ 
in (iii-c). If this is proved, then 
${\mathrm{Exc}}(\tau)$ is locally the product of a 
one-dimensional disk $\Delta^1$ and a 
tree of ${\bold P}^1$'s. 

Therefore, 
$\Lambda$ is smooth at $p_0 \in \Lambda_0$ 
\footnote{The proof is as follows. Let $\Sigma a_iE_i$ 
be the fundamental cycle in the sense of Artin. 
Let $i_1$ be an index which attain the maximal 
coefficient in the cycle. Consider a 
$\bold Q$-divisor $G := \Sigma (a_i/a_{i_1})E_i$. 
Define $\lceil -G \rceil := 
\Sigma \lceil -a_i/a_{i_1} \rceil E_i$, 
where $\lceil -a_i/a_{i_1} \rceil$ denote 
the smallest integer $r$ satisfying 
$r \geq -a_i/a_{i_1}$.
By definition, we have 
$\lceil -G \rceil  = 
-\Sigma_{1 \leq l \leq k} E_{i_l}$, 
where $i_1$, $i_2$, ..., $i_k$ run through all indices 
which attain the maximal coefficient in the cycle. 
Since $-G$ is $\tau$-nef, 
$R^1\tau_*{\cal O}_{\tilde  H}(\lceil -G \rceil) = 0$ 
by Kawamata-Viehweg vanishing theorem. 
For simplicity, we write 
$E := \Sigma_{1 \leq i \leq m}E_i$ and  
$E' := E + \lceil -G \rceil$. 
By the exact sequence 
$\tau_*{\cal O}_{E'}(\lceil -G \rceil) \to 
R^1\tau_*{\cal O}_{\tilde H}(-E) \to 
R^1\tau_*{\cal O}_{\tilde H}(\lceil -G \rceil)$, 
we see that $R^1\tau_*{\cal O}_{\tilde H}(-E) = 0$ 
because the 1-st term vanishes in the sequence. 
The natural map $\tau_*{\cal O}_{\tilde H} 
\to \tau_*{\cal O}_{E}$ is therefore a surjection, 
which implies that $\Lambda$ is a 
normal curve, hence is smooth.}.

Assume that $E_i \cap E_j$ is a section of 
$p_i : E_i \to \Delta^1$, but is a multiple 
section of $p_j : E_j \to \Delta^1$ of degree $> 1$. 
Let $Q \in C_i \cap C_j$. 
Since $N_{C_j/\tilde X} \cong 
{\cal O}^{\oplus (n-2)}\oplus{\cal O}(-2)$ 
and since ${\mathrm{Hilb}}(\tilde X)$ 
is smooth at [$C_j$], 
we can take an (n-2) dimensional subvariety $D$ 
of $\tilde X$ (at least locally around $C_j$) 
such that (1): $D$ has a ${\bold P}^1$-bundle 
structure over an (n-3)-dimensional disc 
$\Delta^{n-3}$ which contains $C_j$ as a central fiber, 
(2): $D$ meets $\tilde H$ in $C_j$ normally. 
Take (n-3) coordinate axis $\Delta_k$ 
($1 \leq k \leq n-3$) of $\Delta^{n-3}$. 
By pulling back the ${\bold P}^1$-bundle 
$D \to \Delta^{n-3}$ by $\Delta_k \to \Delta^{n-3}$, 
we have a smooth surface $D_k$ which has a 
${\bold P}^1$-bundle structure over $\Delta_k$.

We take local coordinates 
$(x,y,z_1,z_2, ..., z_{n-3},t)$ at 
$Q \in \tilde X$ in such a way that \vspace{0.15cm}

(a) ${\bold C}<\partial/\partial x> = 
T_{C_i, Q}$; \vspace{0.13cm}

(b) ${\bold C}<\partial/\partial y> = 
T_{C_j, Q}$; \vspace{0.13cm}

(c) ${\bold C}<\partial/\partial x, 
\partial/\partial y> = T_{E_i, Q}$; \vspace{0.13cm}

(d) ${\bold C}<\partial/\partial y, 
\partial/\partial t> = T_{E_j, Q}$; \vspace{0.13cm}

(e) ${\bold C}<\partial/\partial y, 
\partial/\partial z_k> = T_{D_k, Q}$, 
k = 1,2, ..., n-3. \vspace{0.12cm}
   
 Let $\omega$ be a non-degenerate 2-form on 
$\tilde X$. At $Q$, $\omega_Q$ can be written 
as a linear combination of $dx \wedge dy$, 
$dx \wedge dz_k$, $dx \wedge  dt$, $dy \wedge dz_k$, 
$dy \wedge dt$, $dz_k \wedge dz_l$ and 
$dz_k \wedge dt$.

 Since $E_i$, $E_j$ and $D_k$ ($1 \leq k \leq n-3$) 
are all ${\bold P}^1$-bundles over smooth curves, 
$\omega\vert_{E_i} = \omega\vert_{E_j} = 
\omega\vert_{D_k} = 0$. In particular, 
the terms $dx \wedge dy$, $dy \wedge dz_k$, 
and $dy \wedge dt$ never appear in $\omega_Q$.  
By definition, $\wedge^{n/2} \omega_Q \neq 0$, 
but this is a contradiction.    \vspace{0.15cm}    
  
{\bf (1.6.3)}. {\em there is a non-degenerate 
2-form on $\Sigma\setminus\Sigma_0$}.  \vspace{0.13cm}

 By (1.6.2), 
$\Sigma\setminus\Sigma_0$ is a smooth 
(n-2)-dimensional subvariety. Let 
$E^0 := \pi^{-1}(\Sigma\setminus\Sigma_0)$. 
 $E^0$ is a ${\bold P}^1$-tree 
bundle over $\Sigma\setminus\Sigma_0$. 
The non-degenerate 2-form $\omega$ on $\tilde X$ 
is restricted to a non-zero 2-form $\omega'$ 
on $E^0$. The $\omega'$ must be the pull back of a 
2-form on $\Sigma\setminus\Sigma_0$. 
Since $\wedge^{n/2 - 1}\omega'$ does not vanish 
on the smooth part of $E^0$ 
(cf. Proof of (1.4)), this 2-form on 
$\Sigma\setminus\Sigma_0$ should be 
non-degenerate.  \vspace{0.2cm}

{\bf Examples (1.7)}. {\bf (i)} Let $S$ be a projective 
K3 surface containing a $(-2)$-curve $C$. 
Let $S \to \overline S$ be the birational contraction 
map sending $C$ to a point $p \in \overline S$. 
Let $\tilde X := {\mathrm{Hilb}}^2(S)$ 
be the Hilbert scheme 
parametrizing lenth 2 points on $S$. 
Note that $\tilde X$ is a symplectic 4-fold 
obtained as a resolution of the symmetric product 
${\mathrm{Sym}}^2(S) := S \times S / {\bold Z}_2$ of 
$S$ (cf. [Fu]). Let $X$ be the symmetric 
product ${\mathrm{Sym}}^2(\overline S)$ of $\overline S$.  
Then there is a birational projective morphism 
$\pi : \tilde X \to X$.  The singular locus 
$\Sigma$ of $X$ consists of two irreducible 
components $\Sigma^{(1)}$ and $\Sigma^{(2)}$, 
where both of them are isomorphic to 
$\overline S$ and 
$\Sigma^{(1)} \cap \Sigma^{(2)} = \{(p,p)\}$. 
We can take the 0-dimensional subvariety 
$\{(p,p)\}$ as the $\Sigma_0$ in (1.6). 
$X$ has $A_1$ singularities along 
$\Sigma$ except at $\{(p,p)\}$. 
The fiber $\pi^{-1}((p,p))$ has two irreducible 
components which are isomorphic to 
${\bold P}^2$ and the Hirzebruch surface 
$F_1$. \vspace{0.15cm}

{\bf (ii)} This is an example of a small birational 
contraction map $\pi : \tilde X \to X$ where 
its flop ( = elementary transformation) 
does not preserve the projectivity 
(and also Kaehlerity). In particular, 
$\pi$ is not a projective morphism.

Let $S \to {\bold P}^1$ be an elliptic, 
projective K3 surface which has two type-$I_3$ 
singular fibers 
(cycle of three smooth rational curves). 
Denote by $C = C_1 + C_2 + C_3$ and by 
$D_1 + D_2 + D_3$ these singular fibers. 
Let $\tilde X := {\mathrm{Hilb}}^2(S)$ be the 
Hilbert scheme parametrizing length 2 points on $S$. 
There is a birational projective morphism 
$\mu: \tilde X \to {\mathrm{Sym}}^2(S)$.  
Let $F_i$ (resp. $G_i$) be the proper transform 
of ${\mathrm{Sym}}^2(C_i) \subset {\mathrm{Sym}}^2(S)$ 
(resp. ${\mathrm{Sym}}^2(D_i) 
\subset {\mathrm{Sym}}^2(S)$) by $\mu$.  
The $F_i$'s and $G_i$'s are mutually disjoint and 
they are isomorphic to $\bold P^2$. 
Sinse $N_{F_i/\tilde X} \cong \Omega^1_{\bold P^2}$, 
there exists a (not necessarily projective) 
birational map $\pi : \tilde X \to X$ 
which contracts $F_1$, $F_2$ and $F_3$ to points. 
Let $\pi^+ : \tilde X^+ \to X$ be the flop 
of $\pi$. $\tilde X^+$ is 
obtained by elementary transformations 
along $F_i$'s. 

We shall prove that $\tilde X^+$ is non-projective 
(hence is non-Kaehler because $\tilde X^+$ 
is a Moishezon variety). Let $l_i$ be 
a line on  $F_i$  and $m_i$ a line on $G_i$. 
Then $l_1 + l_2 + l_3$ 
is numerically equivalent to $m_1 + m_2 + m_3$ 
as an algebraic 1-cycle 
on $\tilde X$. In fact, let 
$L \in \mathrm{Pic}(\tilde X)$. 
Then $L \cong {\cal O}_{\tilde X}(aE) 
\otimes \mu^*M$, where $E$ is the 
$\mu$-exceptional divisor, 
$M \in \mathrm{Pic}({\mathrm{Sym}}^2(S))$ and $a \in 
\bold Z$. Let $p: S \times S \to {\mathrm{Sym}}^2(S)$ 
be the natural Galois cover 
with Galois group $G = \bold Z/2\bold Z$. 
Let $q_j : S \times S \to S$ 
be the j-th projection $(j = 1, 2)$.  
Since $H^1(S, {\cal O}_S) = 0$, we 
can write $p^*M = q_1^*M_1 
\otimes q_2^*M_2$ with $M_1$, $M_2 \in 
\mathrm{Pic}(S)$. Since $p^*M$ is G-invariant, 
$M_1 = M_2$. Therefore

$$  (L.l_i)_{\tilde X} = 
(aE. l_i)_{\tilde X} + (M_1. C_i)_S = 
a + (M_1. C_i)_S $$ 

$$  (L.m_i)_{\tilde X} = 
(aE. m_i)_{\tilde X} + (M_1. D_i)_S = 
a + (M_1. D_i)_S $$

As $C_1 + C_2 + C_3$ is linearly equivalent 
to $D_1 + D_2 + D_3$ on $S$, 
it follows that $(L. l_1 + l_2 + l_3)_{\tilde X} 
= (L. m_1 + m_2 + m_3)_{\tilde X}$. 
Therefore $l_1 + l_2 + l_3$ is numerically 
equivalent to $m_1 + m_2 + m_3$. 
We shall derive a contradiction 
assuming that $\tilde X^+$ is projective. 
Let $H$ be an ample divisor 
of $\tilde X^+$. Let 
$H' \in \mathrm{Pic}(\tilde X)$ be the strict 
transform of $H$ by the birational map 
$\phi : \tilde X^+ - - \to \tilde X$. 
By definition of an elementary transformation 
$H'$ is negative along each $F_i$. 
$H'$ is positive along each $G_i$ because 
$\phi$ is an isomorphism around each $G_i$. 
In particular, 
$(H'. l_i) < 0$, and $(H'. m_i) > 0$. 
Hence $(H'. l_1 + l_2 + l_3) < 0$ and 
$(H'. m_1 + m_2 + m_3) > 0$, which is a contradiction 
because $l_1 + l_2 + l_3$ is numerically equivalent 
to $m_1 + m_2 + m_3$. 
Finally note that $\pi : \tilde X \to X$ 
is not a projective morphism 
because if so, then $X$ is projective and 
hence $\tilde X^+ = 
{\mathrm{Proj}}_X(\oplus_m 
\pi_*{\cal O}_{\tilde X}(mH'))$ 
is also projective, which is a contradiction. 
\vspace{0.2cm}          

{\bf (iii)}  Let $A$ be an Abelian surface, and let 
${\mathrm{Hilb}}^3(A)$ be the Hilbert scheme parametrizing 
length 3 points on $A$. The Albanese map 
$Alb: {\mathrm{Hilb}}^3(A) \to A$ factors through the 
symmetric product 
${\mathrm{Sym}}^3(A) := A \times A \times A / S_3$ 
as ${\mathrm{Hilb}}^3(A) \to 
{\mathrm{Sym}}^3(A) \stackrel{f}\to A$. 
For $(x,y,z) \in {\mathrm{Sym}}^3(A)$, 
$f(x,y,z) = x + y + z \in A$. 
Take the origin $0 \in A$, and set 
$\tilde X := {\mathrm{Alb}}^{-1}(0)$ and $X := f^{-1}(0)$. 
The $\tilde X$ is called a symplectic manifold 
of Kummer type, and is often denoted by 
${\mathrm{Kum}}^2(A)$. 
There is a birational projective morphism 
$\pi : \tilde X \to X$. Note that $\tilde X$ is a 
symplectic 4-fold (cf. [Be]). The singular locus 
$\Sigma$ of $X$ is isomorphic to $A$, and 
$X$ has $A_1$ singularities along $\Sigma$ 
except at 81 points $\{p_i\}$ ($1 \leq i \leq 81$). 
The fiber $\pi^{-1}(p_i)$ is homeomorphic to 
the normal surface $\overline F_3$ obtained from 
the Hirzebruch surface $F_3$ by contracting 
$(-3)$-curve to a point (cf. [Br]). 
We can take these 81 points as the $\Sigma_0$ 
in (1.6)  \vspace{0.15cm}  

{\bf (iv)} Let $V$ be a smooth cubic 4-fold in 
$\bold P^5$. Let Y be the Hilbert scheme parametrizing 
lines contained in $V$, which is called classically 
a Fano scheme. Then $Y$ is a symplectic manifold 
of dimension 4. Moreover, $Y$ is deformation equivalent 
to ${\mathrm{Hilb}}^2(S)$ which parametrizes length 2 
points of a K3 surface $S$ ([B-D]). 

We choose a cubic 4-fold $V$ defined by the equation 
$f(T_0, T_2, T_4) + g(T_1, T_3, T_5) = 0$ where $T_i$'s 
are homogenous coordinates of $\bold P^5$. The 
cyclic group $G = \bold Z/ 3\bold Z$ acts on $\bold P^5$ 
by $T_0 \to T_0$, $T_1 \to \zeta T_1$, $T_2 \to T_2$, 
$T_3 \to \zeta T_3$, $T_4 \to T_4$ and 
$T_5 \to \zeta T_5$, where $\zeta$ is a primitive 3 
root of 1. $G$ acts also on $V$, and hence naturally 
on $Y$. By using a $G$ equivariant isomorphism 
$H^1(V, \Omega^3_V) 
\cong H^0(Y, \Omega^2_Y)$ ([B-D]), we know that this 
$G$ action preserves a symplectic 2-form on $Y$.

Let us observe the $G$ action on $Y$ in more detail.
We denote by $P$ 
and $P'$ the projective planes defined by $T_1 = T_3 
= T_5 = 0$ and $T_0 = T_2 = T_4 = 0$ respectively. 
Define $C$ to be the cubic curve on $P$ defined by 
$f(T_0, T_2, T_4) = 0$ and define $D$ to be the 
cubic curve on $P'$ defined by $g(T_1, T_3, T_5) = 0$. 
The fixed locus $F$ of 
the $G$ action on $Y$ is the set of 
lines which join 
two points $p \in C$ and $q \in D$. Hence 
$F \cong C \times D$.  

We put $X := Y/G$. Then $X$ is a symplectic V-manifold 
of dim 4 and its singular locus $\Sigma$ 
is isomorphic to $F$. $X$ has $A_2$ singularities 
along $\Sigma$. Then we can take a symplectic resolution 
$\pi : \tilde X \to X$. It is checked that $\tilde X$ 
is birationally equivalent to 
${\mathrm{Kum}}^2(C \times D)$. 
\vspace{0.2cm}   

(1.8) {\em The Structure of the Singular Locus}:

Let $\pi: \tilde X \to X$, $\Sigma$ and $\Sigma_0$
 be the same as (1.6).
Set $U := X \setminus \Sigma_0$, 
$\tilde U = \pi^{-1}(U)$, $\pi_U := 
\pi\vert_{\tilde U}$ and $D_U := {\mathrm{Exc}}(\pi_U)$. 

We shall study the structure of the sheaf $T^1_U := 
{\underline{Ext}}^1({\Omega^1}_U, {\cal O}_U)$. Note that 
$T^1_U$ has support on $\Sigma\setminus\Sigma_0$.

\vspace{0.15cm}

 Let $D_1, ..., D_m$ be the irreducible 
components of $D_U$. 
 $D_U$ is a divisor with 
normally crossing double points. 
Write $D_{i,j}$ for $D_i \cap D_j$.

We shall describe the possible configurations of 
$D_i$'s over each 
connected component of 
$\Sigma\setminus\Sigma_0$. 
We shall assume, for simplicity, 
that $\Sigma\setminus\Sigma_0$ is 
connected.

Let $\pi_1 : U_1 \to U$ be the blowing-up 
with reduced center 
$\Sigma\setminus\Sigma_0$. 
Then we have: \vspace{0.12cm} 

($A_1$): If $U$ has $A_1$ singularities along 
$\Sigma\setminus\Sigma_0$, then 
${\mathrm{Exc}}(\pi_1)$ is a ${\bold P}^1$-bundle 
over $\Sigma\setminus\Sigma_0$. 
In this case $\tilde U = U_1$ and $\pi_U = \pi_1$.  
\vspace{0.12cm}

 ($A_n$), ($n \geq 2$): If $U$ has $A_n$ singularities 
along $\Sigma\setminus\Sigma_0$, 
then ${\mathrm{Exc}}(\pi_1)$ has a fibration over 
$\Sigma\setminus\Sigma_0$ whose general fiber is a 
(reducible) pair of two lines. The $U_1$ has 
$A_{n-2}$ singularities along the double points of 
${\mathrm{Exc}}(\pi_1)$. Note that
${\mathrm{Exc}}(\pi_1)$ is possibly {\em irreducible} 
when $\Sigma\setminus\Sigma_0$ has non-trivial 
fundamental group. 
    \vspace{0.12cm} 

($D_4$):  If $U$ has $D_4$ singularities along 
$\Sigma\setminus\Sigma_0$, then 
${\mathrm{Exc}}(\pi_1)$ is a ${\bold P}^1$-bundle 
over $\Sigma\setminus\Sigma_0$. 
There is an etale multiple section 
$S \subset {\mathrm{Exc}}(\pi_1)$ of degree 3 
such that $U_1$ has $A_1$ singularities along $S$. 
The $S$ possibly has one, two or 
three connected components.     \vspace{0.12cm}

($D_n$) ($n \geq 5$): If $U$ has $D_n$ ($n \geq 5$) 
singularities along $\Sigma\setminus\Sigma_0$, 
then ${\mathrm{Exc}}(\pi_1)$ is a ${\bold P}^1$-bundle 
over $\Sigma\setminus\Sigma_0$. 
There are two disjoint sections $S_1$ and $S_2$ such 
that $U_1$ has $A_1$ singularities along $S_1$ and has 
$D_{n-2}$ ($A_3$ when $n = 5$) 
singularities along $S_2$. 
\vspace{0.12cm}

($E_6$): If $U$ has $E_6$ singularities along 
$\Sigma\setminus\Sigma_0$, then ${\mathrm{Exc}}(\pi_1)$ 
is a ${\bold P}^1$-bundle over 
$\Sigma\setminus\Sigma_0$. There is a section 
$S$ along which $U_1$ has $A_5$ singularities. 
\vspace{0.12cm}

($E_7$): If $U$ has $E_7$ singularities along 
$\Sigma\setminus\Sigma_0$, then ${\mathrm{Exc}}(\pi_1)$ 
is a ${\bold P}^1$-bundle over $\Sigma\setminus\Sigma_0$.
There is a section $S$ along which $U_1$ has 
$D_6$ singularities. 
\vspace{0.12cm}

($E_8$): If $U$ has $E_8$ singularities along 
$\Sigma\setminus\Sigma_0$, then ${\mathrm{Exc}}(\pi_1)$ 
is a ${\bold P}^1$-bundle over $\Sigma\setminus\Sigma_0$. 
There is a section $S$ along which $U_1$ has $E_7$ 
singularities. \vspace{0.12cm}

Successive blowing ups with singular locus 
give us the (minimal) resolution 
$\pi_U : \tilde U = U_k \to 
U_{k-1} \to ... \to U_1 \to U$. 
Note that each step is essentially the same 
as one of $\pi_1$'s described above. 
We can explicitly check that the number $m$ of 
irreducible components of $D_U$ are as follows according 
to the type of singularities of $U$. \vspace{0.15cm}

$(A_n) => m = n$ or $n - [n/2]$.

$(D_4) => m = 4, 3$ or $2$.

$(D_n)$ ($n \geq 5$) $=> m = n$ or $n-1$.

$(E_6) => m = 6$ or $4$.

$(E_7) => m = 7$. 

$(E_8) => m = 8$. \vspace{0.12cm}

A pair of the type of singularities of $U$ and 
the number $m$ is called a type 
of $U$. For example, if $U$ has $A_n$ singularities 
along $\Sigma\setminus\Sigma_0$ and $m = n - [n/2]$, 
then $U$ is of type $(A_n, n - [n/2])$.

Let $I$ be the defining ideal of the 
reduced subscheme $\Sigma\setminus\Sigma_0$ of $U$. 
Denote by $\Sigma^{(n)}$ the subscheme of $U$ 
defined by $I^{n+1}$. Then there is a sequence of 
subschemes supported at $\Sigma\setminus\Sigma_0$: 
$\Sigma^{(0)} \subset \Sigma^{(1)} \subset 
\Sigma^{(2)} \subset ...$. \vspace{0.15cm}

{\bf Lemma (1.9)}. {\em The sheaf $T^1_U$ is 
described as a sequence of extensions according 
to the type of $U$ in the following way:} 
\vspace{0.15cm}

$(A_1)$: $T^1_U \cong T^1_U\vert_{\Sigma^{(0)}} 
\cong {\cal O}_{\Sigma^{(0)}}$
  \vspace{0.12cm}

$(A_n)$, ($n \geq 2$): \vspace{0.12cm}

$$ 0 \to L \to T^1_U\vert_{\Sigma^{(1)}} \to 
{\cal O}_{\Sigma^{(0)}} \to 0 $$

$$ 0 \to L^{\otimes 2} \to T^1_U\vert_{\Sigma^{(2)}} 
\to T^1_U\vert_{\Sigma^{(1)}} \to 0 $$

$$ : $$

$$ 0 \to L^{\otimes n-1} \to T^1_U\vert_{\Sigma^{(n-1)}} 
\to T^1_U\vert_{\Sigma^{(n-2)}} \to 0$$

$$ T^1_U \cong T^1_U\vert_{\Sigma^{(n-1)}} $$

{\em where $L$ is non-trivial line bundle on 
$\Sigma^{(0)}$ with $L^{\otimes 2} \cong 
{\cal O}_{\Sigma^{(0)}}$ if $U$ is of type 
$(A_n, n-[n/2])$, and where $L \cong 
{\cal O}_{\Sigma^{(0)}}$ if $U$ is of type 
$(A_n, n)$}. \vspace{0.12cm}

$(D_4)$: \vspace{0.12cm}

$$ 0 \to E \to T^1_U\vert_{\Sigma^{(1)}} \to 
{\cal O}_{\Sigma^{(0)}} \to 0 $$

$$ 0 \to {\cal O}_{\Sigma^{(0)}} \to 
T^1_U\vert_{\Sigma^{(2)}} \to 
T^1_U\vert_{\Sigma^{(1)}} \to 0 $$

$$ T^1_U \cong T^1_U\vert_{\Sigma^{(2)}} $$

{\em where $E$ is a vector bundle of rank 2 with 
$h^0(E) = 0$ if $U$ is of type $(D_4, 2)$, where 
$E \cong {\cal O}_{\Sigma^{(0)}}\oplus L$ 
with a non-trivial line bundle $L$, 
$L^{\otimes 2} \cong {\cal O}_{\Sigma^{(0)}}$ if 
$U$ is of type $(D_4, 3)$, and where 
$E \cong {\cal O}^{\oplus 2}_{\Sigma^{(0)}}$ if 
$U$ is of type $(D_4, 4)$}. \vspace{0.12cm}

$(D_n)$ ($n \geq 5$): \vspace{0.12cm}

$$  0 \to {\cal O}_{\Sigma^{(0)}}\oplus L \to 
T^1_U\vert_{\Sigma^{(1)}} \to 
{\cal O}_{\Sigma^{(0)}} \to 0 $$

$$ 0 \to {\cal O}_{\Sigma^{(0)}} \to 
T^1_U\vert_{\Sigma^{(2)}} \to 
T^1_U\vert_{\Sigma^{(1)}} \to 0 $$

$$ : $$ 
$$ 0 \to {\cal O}_{\Sigma^{(0)}} 
\to T^1_U\vert_{\Sigma^{(n-2)}} \to 
T^1_U\vert_{\Sigma^{(n-3)}} \to 0$$

$$ T^1_U \cong T^1_U\vert_{\Sigma^{(n-2)}} $$ 

{\em where $L$ is a non-trivial line bundle with 
$L^{\otimes 2} \cong {\cal O}_{\Sigma^{(0)}}$ 
if $U$ is of type $(D_n, n-1)$ and where 
$L \cong {\cal O}_{\Sigma^{(0)}}$ if 
$U$ is of type $(D_n, n)$}.  \vspace{0.12cm}

$(E_6)$: \vspace{0.12cm}

$$ 0 \to E_1 \to T^1_U\vert_{\Sigma^{(1)}} 
\to {\cal O}_{\Sigma^{(0)}} \to 0 $$

$$ 0 \to E_2 \to T^1_U\vert_{\Sigma^{(2)}} 
\to T^1_U\vert_{\Sigma^{(1)}} \to 0 $$

$$ 0 \to {\cal O}_{\Sigma^{(0)}} \to 
T^1_U\vert_{\Sigma^{(3)}} \to 
T^1_U\vert_{\Sigma^{(2)}} \to 0$$

$$ T^1_U \cong T^1_U\vert_{\Sigma^{(3)}} $$ 

{\em where $E_1$ and $E_2$ are vector bundles 
on $\Sigma^{(0)}$ of rank 2 obtained as 
the following extensions} \vspace{0.12cm}

$$ 0 \to {\cal O}_{\Sigma^{(0)}} 
\to E_1 \to L \to 0$$

$$ 0 \to L \to E_2 
\to {\cal O}_{\Sigma^{(0)}} \to 0 $$

{\em where $L$ is a non-trivial line bundle 
with $L^{\otimes 2} \cong 
{\cal O}_{\Sigma^{(0)}}$ if $U$ is of type 
$(E_6, 4)$, and where $L \cong 
{\cal O}_{\Sigma^{(0)}}$ if $U$ is of type 
$(E_6, 6)$}. \vspace{0.12cm}

$(E_7, 7)$, $(E_8, 8)$: \vspace{0.12cm}

$$ 0 \to E_1 \to T^1_U\vert_{\Sigma^(1)} \to 
{\cal O}_{\Sigma^{(0)}} \to 0 $$

$$ 0 \to E_2 \to T^1_U\vert_{\Sigma^{(2)}} \to 
T^1_U\vert_{\Sigma^{(1)}} \to 0 $$

$$ 0 \to E_3 \to T^1_U\vert_{\Sigma^{(3)}} 
\to T^1_U\vert_{\Sigma^{(2)}} \to 0 $$

$$ 0 \to {\cal O}_{\Sigma^{(0)}} \to 
T^1_U\vert_{\Sigma^{(4)}} \to 
T^1_U\vert_{\Sigma^{(3)}} \to 0$$

$$ T^1_U \cong T^1_U\vert_{\Sigma^{(4)}} $$ 
\vspace{0.12cm}

{\em where $E_1$ and $E_2$ are vector bundles of 
rank 2 obtained as extensions of trivial line bundles 
: $0 \to {\cal O}_{\Sigma^{(0)}} \to E_i \to 
{\cal O}_{\Sigma^{(0)}} \to 0 $.  
Moreover, if $U$ is of type $(E_7, 7)$, then 
$E_3 \cong {\cal O}_{\Sigma^{(0)}}$, and if 
$U$ is of type $(E_8, 8)$, then $E_3$ is a 
vector bundle of rank 2 obtained as an extension of 
trivial line bundles : $0 \to {\cal O}_{\Sigma^{(0)}} 
\to E_3 \to {\cal O}_{\Sigma^{(0)}} \to 0$.} 

\vspace{0.12cm}

The proof of this lemma is omitted; I will 
write it elsewhere. Note that we need two 
facts to prove it: $\omega_U \cong {\cal O}_U$ 
and $\omega_{\Sigma^{(0)}} 
\cong {\cal O}_{\Sigma^{(0)}}$. 
 An important collorary of (1.9) 
is the following: \vspace{0.15cm}

{\bf Corollary (1.10)}. 
$$h^0(U, T^1_U) \leq m.$$ 

{\em Proof}. Note that, by Proposition (1.6), 
$\Sigma\setminus\Sigma_0$ is compactified to a 
proper normal variety $\overline\Sigma$ by 
adding codimension 2 points. We see that 
$h^0(U, T^1_U) \leq m$ by the extensions in (1.9)
 in any case. 
\vspace{0.2cm}
 
We shall finally state related results to (1.4), 
which will not be used in the later.
Proposition (1.4) is also available 
if we replace $\pi : \tilde X \to X$ by a projective 
symplectic resolution $\phi : \tilde V \to V$ 
of the germ of a rational Gorenstein singularity 
$0 \in V$ of even dimension $n$. 
When $V$ is an isolated singularity, 
we have the following. \vspace{0.2cm}

{\bf Proposition (1.11)}. {\em Let 
$\phi : \tilde V \to V$ be a 
projective symplectic resolution of 
the germ of an isolated rational 
Gorenstein singularity $0 \in V$ of dimension 
$n \geq 4$, that is, 
$\tilde V$ admits a non-degenerate holomorphic 2-form. 
Then every irreducible component of the exceptional 
locus ${\mathrm{Exc}}(\phi)$ is 
Lagrangian.} \vspace{0.2cm}

{\em Proof}.  First we shall prove: \vspace{0.15cm}

{\bf Lemma (1.12)}. {\em Let 
$\phi : \tilde V \to V$ be a projective, 
crepant resolution of an isolated rational 
Gorenstein singularity of 
even dimension $n$. Then every 
irreducible component of $E := 
{\mathrm{Exc}}(\phi)$ has dimension $\geq n/2$.}  
\vspace{0.2cm}

{\em Proof of (1.12)}. Take an effective divisor 
$\Delta$ of $\tilde X$ 
in such a way that $-\Delta$ is $\phi$-ample. 
If $\epsilon > 0$ is a 
sufficiently small rational number, 
then $(\tilde X, \epsilon\Delta)$ 
is log terminal in the sense of [Ka 2]. 
Since $K_{\tilde X} \sim 0$, 
every irreducible component of $E$ 
is covered by a family of rational 
curves by [Ka 2, Theorem 1]. 

We shall now use the terminology in 
[Ko, Chap. IV].
Let $E_i$ be an irreducible component of $E$. 
Let ${\mathrm{Hom}}_{bir}({\bold P}^1, E_i)$ 
be the Hom scheme parametrizing 
the morphisms from ${\bold P}^1$ to $E_i$ 
which are birational onto 
their images. Let 
${\mathrm{Hom}}^n_{bir}({\bold P}^1, E_i)$ 
be the normalization of 
${\mathrm{Hom}}_{bir}({\bold P}^1, E_i)$. 
By [Ko 1, Chap.IV, Theorem 2.4], there is an 
irreducible component $W_i$ of 
${\mathrm{Hom}}^n_{bir}({\bold P}^1, E_i)$ such that 
$W_i$ is a {\em generically unsplit family} of 
morphisms and such that 
$\overline{{\mathrm{Locus}}(W_i)} = E_i$, where 
${\mathrm{Locus}}(W_i)$ is the locus 
where the images of the morphisms in $W_i$ 
sweep out and $\overline{{\mathrm{
Locus}}(W_i)}$ is its closure.
 
Let $[f] \in W_i$ be a general point. 
Then $W_i$ is also an 
irreducible component of $
{\mathrm{Hom}}^n_{bir}({\bold P}^1, \tilde V)$ 
at $[f]$. We know that 
$\dim_{[f]} W_i \geq \chi ({\bold P}^1, 
f^*\Theta_{\tilde V}) = n$. For the generically 
unsplit family $W_i$ 
of morphisms, we can estimate 
${\mathrm{codim}}({\mathrm{Locus}}(W_i) 
\subset \tilde V)$ 
(cf. [Io, Theorem (0.4), Ko, Chap.IV, Proposition 
(2.5)]). The result is 

$$  {\mathrm{codim}}({\mathrm{Locus}}(W_i) 
\subset \tilde V) \leq n/2 + 
1/2. $$  
    
Note here that $n$ is even.
Since $\overline{{\mathrm{Locus}}(W_i)} = E_i$, 
we have $\dim E_i \geq n/2$. \hspace{0.4cm} Q.E.D.  
\vspace{0.2cm}

{\em Proof of Proposition (1.11) continued}. 
By combining Lemma (1.12) with Proposition (1.4) 
we have  $\dim E_i = n/2$ for each irreducible 
component $E_i$ of $E$.  Let us prove that $E_i$ 
are all Lagrangian. Let $\omega$ be a 
non-degenerate 2-form on $\tilde V$. 
Assume that $\omega\vert_{E_i} \neq 0$ for some $E_i$. 
Take a birational projective morphism 
$\nu : Y \to \tilde V$ in such a way that 
$\nu^{-1}(E)$ becomes a divisor of a 
smooth n-fold $Y$ with normal crossings. 
Write $g = \phi\circ\nu$ for short. 
Then we have 
$H^0(g^{-1}(0), \hat\Omega^2_{g^{-1}(0)}) \neq 0$. 
This contradicts Lemma (1.2). \hspace{0.4cm} Q.E.D.   
\vspace{0.2cm} 
    
\begin{center}
{\bf 2. Deformation theory}
\end{center}

We shall review some generalities of deformation 
theory. For a compact complex space $X$ we denote 
by ${\mathrm{Def}}(X)$ the Kuranishi space of $X$. 
By definition, there is a reference point 
$0 \in {\mathrm{Def}}(X)$ and there is a 
semi-universal flat deformation 
$f: \cal X \to {\mathrm{Def}}(X)$ of 
$X$ with $f^{-1}(0) = X$. When $X$ is reduced, 
the tangent space $T_{{\mathrm{Def}}(X), 0}$ is 
canonically isomorphic to 
$\mathrm{Ext}^1(\Omega^1_X, {\cal O}_X)$. We abbreviate 
this space by ${\bold T}^1_X$. 

Let 
$\pi : \tilde X \to X$ be a proper bimeromorphic 
map of compact complex spaces. Assume that 
$R^1\pi_*{\cal O}_{\tilde X} = 0$ and 
$\pi_*{\cal O}_{\tilde X} = {\cal O}_X$. Then there 
is a natural map ${\mathrm{Def}}(\tilde X) \to 
{\mathrm{Def}}(X)$ (cf. [Ko-Mo, 11.4]).  This map 
naturally induces a map $\pi_* : {\bold T}^1_{\tilde X} 
\to {\bold T}^1_X$. Assume that $\tilde X$ and $X$ 
are both reduced. Then $\pi_*$ is obtained as follows. 

Let $0 \to {\cal O}_{\tilde X} 
\to F \to \Omega^1_{\tilde X} \to 0$ be the 
extension corresponding to an element of 
${\bold T}^1_{\tilde X}$. Operate $\pi_*$ on this 
sequence. Then we have an exact sequence 
$0 \to {\cal O}_X \to \pi_*F 
\to \pi_*\Omega^1_{\tilde X} 
\to 0$ because $R^1\pi_*{\cal O}_{\tilde X} = 0$ and 
$\pi_*{\cal O}_{\tilde X} = {\cal O}_X$. This 
extension gives an element of ${\mathrm{Ext}}^1
(\pi_*\Omega^1_{\tilde X}, {\cal O}_X)$.
By the natural 
map $\Omega^1_X \to \pi_*\Omega^1_{\tilde X}$, 
we obtain an element of ${\bold T}^1_X$. \vspace{0.12cm}

In the remainder of this section, $\tilde X$ is a smooth 
projective symplectic n-fold with $n \geq 4$ and 
$\pi : \tilde X \rightarrow X$ is a birational 
projective morphism from $\tilde X$ to a normal 
n-fold $X$. We shall use the same notation as 
(1.6). Set $U := X\setminus\Sigma_0$.  
\vspace{0.15cm}

{\bf Proposition (2.1)}. {\em There is a commutaive 
diagram} 

\begin{equation}
\begin{CD}
H^1(\tilde X, \Theta_{\tilde X}) @>>> 
H^1(\pi^{-1}(U), \Theta_{\pi^{-1}(U)}) \\
@VVV @VVV \\
{\bold T}^1_X @>>> {\bold T}^1_U 
\end{CD}
\end{equation}

{\em and the horizontal maps are both isomorphisms}.  
\vspace{0.15cm}

{\em Proof}. Since $X$ has rational singularities, 
$R^1\pi_*{\cal O}_{\tilde X} = 0$ and hence we have 
vertical maps by [Bi, Wa]. The horizontal map on 
the second row is an isomorphism by the argument of 
[Na 4, Propositions (1.1), (1.2)] or [Ko-Mo, (12.5.6)]. 
We only have to prove that the horizontal map 
on the first row is an isomorphism. 

Since we will treat non-compact varieties, we 
shall fix the notation here. Assume that a complex 
space $W$ admits a structure of an algebraic 
scheme $W^{alg}$ over $\bold C$. Let $F$ be a coherent 
analytic sheaf $F$ on $W$ which comes 
from an algebraic 
coherent sheaf $F^{alg}$ on $W^{alg}$. 
Then we write $H^{*}(W, F^{alg})$ for 
$H^{*}(W^{alg}, F^{alg})$ by abuse of 
notation. Note that there is a natural map 
$H^{*}(W, F^{alg}) \to H^{*}(W, F)$. 
  
We can see that the map $H^1(\pi^{-1}(U), 
\Theta^{alg}_{\pi^{-1}(U)}) \to 
H^1(\pi^{-1}(U), \Theta_{\pi^{-1}(U)})$ is an 
isomorphism (cf. Appendix). Thus,  
let us consider the exact sequence of local cohomology 
in the algebraic category.

$$ H^1_{\pi^{-1}(\Sigma_0)}(\Theta^{alg}_{\tilde X}) \to 
H^1(\Theta^{alg}_{\tilde X}) \to H^1(\pi^{-1}(U), 
\Theta^{alg}_{\pi^{-1}(U)}) \to 
H^2_{\pi^{-1}(\Sigma_0)}(\Theta^{alg}_{\tilde X}). $$

We only have to prove that the middle map is an 
isomorphism by GAGA and by the fact mentioned above. 
Let $X_{\Sigma_0}$ (resp. ${\tilde X}_{\Sigma_0}$) 
be the formal completion of $X$ (resp. $\tilde X$) 
along $\Sigma_0$ (resp. ${\pi^{-1}(\Sigma_0)}$). 
By duality, we have $H^{n-1}({\tilde X}_{\Sigma_0}, 
\Omega^{1, alg}_{{\tilde X}_{\Sigma_0}}) \cong 
H^1_{\pi^{-1}(\Sigma_0)}(\Theta^{alg}_{\tilde X})^*$ and 
$H^{n-2}({\tilde X}_{\Sigma_0}, 
\Omega^{1,alg}_{{\tilde X}_{\Sigma_0}}) \cong 
H^2_{\pi^{-1}(\Sigma_0)}(\Theta^{alg}_{\tilde X})^*$.
     
Note that $\Omega^{1, alg}_{\tilde X} \cong 
\Omega^{n-1, alg}_{\tilde X}$ by a non-degenerate 2-form 
$\omega$. Therefore, we have to prove that 
$H^{n-1}(\tilde X_{\Sigma_0}, 
\Omega^{n-1, alg}_{\tilde X_{\Sigma_0}}) = 
H^{n-2}(\tilde X_{\Sigma_0}, 
\Omega^{n-1, alg}_{\tilde X_{\Sigma_0}}) = 0$.  
We shall prove that 
$R^k{\pi}_*\Omega^{n-1, alg}_{\tilde X} = 0$ if 
$k \geq 2$. If these are proved, then by the Leray 
spectral sequence 

$$  E^{p,q}_2 = H^p(X_{\Sigma_0}, 
R^q{\pi}_*\Omega^{n-1, alg}_{\tilde X}) => 
H^{p+q}(\tilde X_{\Sigma_0}, 
\Omega^{n-1, alg}_{\tilde X_{\Sigma_0}}) $$

we have $H^{n-1}(\tilde X_{\Sigma_0}, 
\Omega^{n-1, alg}_{\tilde X_{\Sigma_0}}) = 
H^{n-2}(\tilde X_{\Sigma_0}, 
\Omega^{n-1, alg}_{\tilde X_{\Sigma_0}}) = 0$ 
because $\dim \Sigma_0 \leq n-4$ by Proposition (1.6). 

By GAGA, it suffices to show, in the 
analytic category, that 
$R^k{\pi}_*\Omega^{n-1}_{\tilde X} = 0$ if 
$k \geq 2$. 
 
Let $\nu : Y \to \tilde X$ be a composition of 
blowing ups with smooth centers such that 
the total transform of ${\mathrm{Exc}}(\pi)$ is a 
divisor with normal crossings. 
Set $f := \pi\circ\nu$. We put $E := {\mathrm{Exc}}(f)$ 
and $E' := {\mathrm{Exc}}(\nu)$. \vspace{0.2cm}

{\bf Claim 1}. $R^k{\pi}_*\Omega^{n-1}_{\tilde X} 
\cong R^kf_*\Omega^{n-1}_Y(log E')(-E')$. 
\vspace{0.15cm}

{\em Proof}.  By [St 2], 
$R^l{\nu}_*\Omega^{n-1}_Y(log E')(-E') = 0$ for 
$l \geq 2$. The same statement also holds 
when $l = 1$. The proof of this fact is the 
following. 
Since ${\nu}_*{\hat \Omega}^{n-1}_{E'} = 0$, 
we have an exact sequence 

$$ 0 \to R^1{\nu}_*\Omega^{n-1}_Y(log E')(-E') \to 
R^1{\nu}_*\Omega^{n-1}_Y 
\stackrel{{\gamma}'}\to 
R^1{\nu}_*{\hat \Omega}^{n-1}_{E'}. $$ 

The map ${\gamma}'$ is an isomorphism because 
$\tilde X$ is smooth and $\nu$ is a composition 
of the blowing-ups with smooth centers. 
In fact, at first, by the exact sequence 

$$ R^1\nu_*\Omega^{n-1}_Y \to 
R^1\nu_*\hat{\Omega}^{n-1}_{E'} \to 
R^2\nu_*\Omega^{n-1}_Y(log E')(-E') $$

${\gamma}'$ is surjective since the last 
term vanishes by [St 2].
Assume that $\nu$ is a composition of exactly 
k blowing-ups. 
We shall prove that ${\gamma}'$ is an isomorphism 
by the induction on k.  
We can check it directly when $k = 1$. 
Assume that $k >1$. 
Decompose $\nu$ as $Y \stackrel{\nu_2}\to Y_1 
\stackrel{\nu_1}\to \tilde X$ in such a way that 
$\nu_1$ is a blowing-up with a smooth center. 
Let $E_1 := {\mathrm{Exc}}(\nu_1)$. 
Then the proper transform $E'_1$ of $E_1$ by 
$\nu_2$ is an irreducible component of 
$E' := {\mathrm{Exc}}(\nu)$. Let $E' = \Sigma {E'}_i$ 
be the irreducible decomposition. 
We see that ${\nu_2}_*{\hat{\Omega}^{n-1}_{E'}} 
\cong \oplus{\nu_2}_*{\Omega}^{n-1}_{{E'}_i} 
\cong {\nu_2}_*{\Omega}^{n-1}_{E'_1} 
\cong {\Omega}^{n-1}_{E_1}$.  
There is a commutative diagram with exact rows 

\begin{equation}
\begin{CD}
0 @>>> R^1{\nu_1}_*{\Omega}^{n-1}_{Y_1} @>>> 
R^1\nu_*\Omega^{n-1}_Y @>>> 
{\nu_1}_*R^1{\nu_2}_*{\Omega^{n-1}_Y} \\
@. @VVV @V{{\gamma}'}VV @VVV \\
0 @>>> R^1{\nu_1}_*{\hat\Omega}^{n-1}_{E_1} 
@>>> R^1\nu_*\hat{\Omega}^{n-1}_{E'} 
@>>> {\nu_1}_*R^1{\nu_2}_*{\hat\Omega}^{n-1}_{E'}
\end{CD}
\end{equation}   
   
The vertical maps except ${\gamma}'$ are both 
isomorphisms by the induction, hence 
${\gamma}'$ is injective by the diagram. 
Thus ${\gamma}'$ is an isomorphism. \vspace{0.15cm}

Now by a Leray spectral sequence, 
$R^kf_*\Omega^{n-1}_Y(log E')(-E') \cong 
R^k{\pi}_*({\nu}_*\Omega^{n-1}_Y(log E')(-E'))$.  
By the exact sequence 

$$ 0 \to \nu_*\Omega^{n-1}_Y(log E')(-E') 
\to \nu_*\Omega^{n-1}_Y \to 
\nu_*{\hat \Omega}^{n-1}_{E'} = 0 $$

we see that ${\nu}_*\Omega^{n-1}_Y(log E')(-E') 
\cong  \nu_*\Omega^{n-1}_Y$, where the second term 
is isomorphic to 
$\Omega^{n-1}_{\tilde X}$ (cf. [St 1]).  Q.E.D. 
\vspace{0.2cm}

{\bf Claim 2}. 
{\em $R^kf_*\Omega^{n-1}_Y(log E')(-E') = 0$ if 
$k \geq 2$}. 
\vspace{0.15cm}

{\em Proof}. Consider the exact sequence 

$$  R^{k-1}f_*\Omega^{n-1}_Y \stackrel{\alpha}\to 
R^{k-1}f_*{\hat \Omega}^{n-1}_{E'} \to 
R^kf_*\Omega^{n-1}_Y(log E')(-E') \to 
R^kf_*\Omega^{n-1}_Y \stackrel{\beta}\to 
R^kf_*{\hat \Omega}^{n-1}_{E'}. $$

The map $\alpha$ is factorized as 
$R^{k-1}f_*\Omega^{n-1}_Y \to 
R^{k-1}f_*{\hat \Omega}^{n-1}_{E} \to 
R^{k-1}f_*{\hat \Omega}^{n-1}_{E'}$, 
and the second map is clearly a surjection. 
The first map is also surjective 
by the exact sequence 

$$ R^{k-1}f_*\Omega^{n-1}_Y \to 
R^{k-1}f_*{\hat \Omega}^{n-1}_{E} 
\to R^kf_*\Omega^{n-1}_Y(log E)(-E) $$ 

because $R^kf_*\Omega^{n-1}_Y(log E)(-E) = 0$ 
for $k \geq 2$ by [St 2]. Hence 
$\alpha$ is a surjection. 

The map $\beta$ is similarly factorized as 
$R^{k}f_*\Omega^{n-1}_Y \to 
R^{k}f_*{\hat \Omega}^{n-1}_{E} \to 
R^{k}f_*{\hat \Omega}^{n-1}_{E'}$. 
Note that the second map is an isomorphism. 
Indeed, let $E''$ be an $f$-exceptional divisor 
which is not contained in $E'$. 
Then, by Corollary (1.5), 
$E''$ is mapped to  an (n-2)-dimensional subvariety 
of $X$ by $f$; in particular, a general fiber 
of $E'' \to f(E'')$ has dimension 1. 
By [Ko 2] $R^kf_*{\hat \Omega}^{n-1}_{E''} = 0$ 
if $k \geq 2$. Hence 
$R^{k}f_*{\hat \Omega}^{n-1}_{E} \cong 
R^{k}f_*{\hat \Omega}^{n-1}_{E'}$.  

The first map is injective by the exact sequence    

$$ R^kf_*\Omega^{n-1}_Y(log E)(-E) \to 
R^{k}f_*\Omega^{n-1}_Y \to 
R^{k}f_*{\hat \Omega}^{n-1}_{E} $$     

because the first term vanishes by [St 2]. 
Hence $\beta$ is an injection. Q.E.D. \vspace{0.2cm}

{\bf Theorem (2.2)}. {\em Let $\pi : \tilde X \to X$ 
be a symplectic resolution of a projective symplectic 
variety $X$ of dimension $n$. Then the Kuranishi spaces 
${\mathrm{Def}}(\tilde X)$ and ${\mathrm{Def}}(X)$ 
are both smooth of the same dimension. 
The natural map $\pi_* : {\mathrm{Def}}(\tilde X) 
\to {\mathrm{Def}}(X)$ is a finite covering. 
Moreover, $X$ has a flat deformation to a 
symplectic n-fold $X_t$. Any smoothing $X_t$ of 
$X$ is obtained as a flat deformation of 
$\tilde X$.} \vspace{0.15cm}

{\em Proof}. By Bogomolov [Bo] the Kuranishi space
${\mathrm{Def}}(\tilde X)$ is smooth for a symplectic 
manifold $\tilde X$. \vspace{0.15cm}

{\bf Claim 1}. $\dim {\bold T}^1_X \leq
\dim {\mathrm{Def}}(\tilde X).$ \vspace{0.12cm}

{\em Proof}. Since ${\mathrm{Def}}(\tilde X)$ is smooth, 
we have to prove that $\dim {\bold T}^1_X
 \leq h^1(\tilde X, \Theta_{\tilde X})$. 
By Proposition (2.1) it suffices to prove that 
$\dim {\bold T}^1_U \leq h^1(\pi^{-1}(U), 
\Theta_{\pi^{-1}(U)})$. We shall use the same 
notation as (1.8). Recall that $\pi_U := 
\pi\vert_{\tilde U}$, $D_U := {\mathrm{Exc}}(\pi_U)$ 
and $D_U = D_1 \cup ... \cup D_m$ is the irreducible 
decomposition. 
\vspace{0.12cm}

(i) By Corollary (1.10) we have $h^0(U, T^1_U) \leq m$. 
\vspace{0.12cm}

(ii) Consider the commutative diagram with exact rows

\begin{equation}
\begin{CD}
0 @>>> H^1(U, \Theta_U) @>>> H^1(\tilde U, 
\Theta_{\tilde U}) @>{\tilde\eta}>> 
H^0(U, R^1(\pi_U)_*\Theta_{\tilde U}) \\
@. @V{id}VV @VVV @VVV \\
0 @>>> H^1(U, \Theta_U) @>>> {\bold T}^1_U 
@>{\eta}>> H^0(U, T^1_U)
\end{CD}
\end{equation}

We shall prove that $\tilde{\eta}$ is 
surjective and prove
that $h^0(U, R^1(\pi_U)_*\Theta_{\tilde U}) = m$,
where $m$ is the number of irreducible 
components of $D_U$.

By a non-degenerate 2-form $\omega$ on $\tilde X$, 
we have $H^0(U,
R^1(\pi_U)_*\Theta_{\tilde U}) \cong 
H^0(U, R^1(\pi_U)_*\Omega^1_{\tilde U})$.

By the definition of $U$ we know that 
$R^1(\pi_U)_*\Omega^1_{\tilde
U} \cong  R^1(\pi_U)_*{\hat\Omega^1_{D_U}}$.
There is an exact sequence

$$ \oplus_i(\pi_U)_*\Omega^1_{\tilde D_i} \to
\oplus_{i>j}(\pi_U)_*\Omega^1_{\tilde D_{i,j}} \to
R^1(\pi_U)_*\hat\Omega^1_{D_U} \to
\oplus_iR^1(\pi_U)_*\Omega^1_{\tilde D_i}$$

where $\tilde D_i$ and $\tilde D_{i,j}$ are 
normalizations of $D_i$ and $D_{i,j}$ respectively. 

The first map is surjective by the description of
$D_U$ in (i).
The last map is a surjection 
(and hence an isomorphism) because
$R^1(\pi_U)_*\Omega^1_{\tilde D_{i,j}}
 = 0$.
Let $\tilde D_i \stackrel{\pi_i}\to S_i 
\stackrel{\varphi_i}\to \Sigma\setminus\Sigma_0$ 
be the Stein factorization. 
Since $\varphi_i$ is finite, 
$R^1(\pi_U)_*\Omega^1_{\tilde D_i} \cong 
(\varphi_i)_*R^1(\pi_i)_*\Omega^1_{\tilde D_i}$. 
Since $\tilde D_i$ is a ${\bold P}^1$-bundle over 
$S_i$, we have $R^1(\pi_i)_*\Omega^1_{\tilde D_i} 
\cong {\cal O}_{S_i}$. Thus 
$H^0(U, R^1(\pi_U)_*\Omega^1_{\tilde D_i}) 
\cong H^0(U, {\phi_i}_*{\cal O}_{S_i}) = 
H^0(S_i, {\cal O}_{S_i}) = \bold{C}$. 

As a consequence, 
$h^0(U, R^1(\pi_U)_*{\Theta_{\tilde U}}) = m$. 
\vspace{0.12cm}

(iii) By a non-degenerate 2-form $\omega$ 
on $\tilde X$, the map $\tilde\eta$ is identified 
with the map $H^1(\tilde U, \Omega^1_{\tilde U}) 
\to H^0(U, R^1(\pi_U)_*\Omega^1_{\tilde U})$. 
We have a commutative diagram 
 
\begin{equation}
\begin{CD}
H^1(\tilde U, {\cal O}^*_{\tilde U}) @>>> 
H^0(U, R^1(\pi_U)_*{\cal O}^*_{\tilde U}) \\
@V{dlog}VV @V{dlog}VV \\
H^1(\tilde U, \Omega^1_{\tilde U}) @>>> 
H^0(U, R^1(\pi_U)_*\Omega^1_{\tilde U})
\end{CD}
\end{equation}

The ${\bold C}$-vector space 
$H^0(U, R^1(\pi_U)_*\Omega^1_{\tilde U})$
is generated by the images of $[D_i] \in 
H^1(\tilde U, {\cal O}^*_{\tilde U})$ 
($1 \leq i \leq m$) by (ii). 
Thus the horizontal map at the bottom 
is surjective, and $\tilde\eta$ is 
also surjective.  \vspace{0.12cm}

(iv) By the commutative diagram (6) 
$h^0(\tilde U, \Theta_{\tilde U}) = 
h^0(U, \Theta_U) + m$ because 
$h^0(U, R^1(\pi_U)_*{\Theta_{\tilde U}}) = m$ 
and $\tilde\eta$ is surjective. 
On the other hand, 
$dim {\bold T}^1_U \leq  
h^0(U, \Theta_U) + m$ because $h^0(U, T^1_U) \leq m$. 
These imply that 
$\dim {\bold T}^1_U \leq 
h^1(\tilde U, \Theta_{\tilde U})$. \vspace{0.2cm}

{\bf Claim 2}. {\em The map $\pi_* : 
{\mathrm{Def}}(\tilde X) \to {\mathrm{Def}}(X)$ 
has finite fibers.}  \vspace{0.12cm}

{\em Proof}. Let 
$\Pi : \tilde{\cal X} \to X \times {\Delta}^1$ 
be a flat deformation of the map 
$\pi : \tilde X \to X$ over a 1-dim disc $\Delta^1$. 
We have to show that $\tilde{\cal X} \cong 
{\tilde X}\times{\Delta}^1$. 
Let $S_n := {\mathrm{Spec}}{\bold C}[t]/(t^{n+1})$ 
and let ${\tilde X}_n$ be the pull back of 
$\tilde{\cal X}$ by the natural embedding 
$S_n \to \Delta^1$. We have to prove that 
${\tilde X}_n \cong {\tilde X}\times S_n$ 
for all $n$. By Proposition (2.1), 
there is a one to one corrspondence between 
infinitesimal deformations of $\tilde X$ and 
infinitesimal deformations of 
$\tilde U := \pi^{-1}(U)$. 
Therefore, it suffices to prove the same 
statement by replacing $X$ by $U$ and 
$\tilde X$ by $\tilde U$ respectively. 
But, then ${\tilde U}_n$ should be a 
(relatively) minimal resolution of 
$U \times S_n$. By the uniqueness of 
minimal resolution, we have 
${\tilde U}_n \cong {\tilde U \times S_n}$. 
\hspace{0.5cm} Q.E.D. \vspace{0.2cm}

By Claims 1 and 2, $\dim {\mathrm{Def}}(\tilde X) = 
\dim {\bold T}^1_X$. Since 
${\mathrm{Def}}(\tilde X)$ is smooth by [Bo], 
${\mathrm{Def}}(X)$ is also smooth. 
Moreover, $\pi_*$ is a finite covering 
(cf. [Fi, 3.2, p 132]). \vspace{0.2cm}

{\bf Claim 3}. {\em $X$ has a flat deformation 
to a smooth symplectic n-fold $X_t$ 
such that $X_t$ is a small deformation of 
$\tilde X$.}  \vspace{0.15cm}

{\em Proof}. The proof is due to [Fu, (3) 
in the proof of Theorem (5.7)]. By the existence of a 
non-degenerate 2-form $\omega$, 
there is an obstruction to extending a 
holomorphic curve on $\tilde X$ sideways 
in a given one-parameter small deformation 
$\tilde{\cal X} \to \Delta^1$. 
Therefore, if we take a general curve of 
${\mathrm{Def}}(\tilde X)$ passing through 
the origin and take a corresponding small 
deformation of $\tilde X$, 
then no holomorphic curves survive. 
A detailed argument on this fact can be found 
in [Fu, Theorem (4.8), (1)]; 
the theorem assumes that $\tilde X$ 
is primitively symplectic, however, 
one can prove the same result in a general 
case by a minor modification.  

Let $t \in {\mathrm{Def}}(X)$ be a generic point 
(that is, $t$ is outside the union of a 
countable number of proper subvarieties of 
${\mathrm{Def}}(X)$).  Since 
$\pi_* : {\mathrm{Def}}(\tilde X) \to {\mathrm{Def}}(X)$ 
is a finite covering, we may assume that 
$X_t$ has a symplectic resolution 
$\pi_t : {\tilde X}_t \to X_t$. 
By the argument above, ${\tilde X}_t$ 
contains no curves. 
By Chow lemma [Hi], there is a bimeromorphic 
projective map $h : W \to X_t$ such that 
$h$ is factored through $\pi_t$. 
Since $h^{-1}(p)$ is the union of 
projective varieties for $p \in X_t$, 
${\pi_t}^{-1}(p)$ is the union of 
Moishezon varieties. 
If $\pi_t$ is not an isomorphism, then 
${\pi_t}^{-1}(p)$ has positive dimension 
for some point $p \in X_t$; hence 
${\tilde X}_t$ contains curves, 
which is a contradiction. Thus $\pi_t$ 
is an isomorphism and $X_t$ is a 
(smooth) symplectic n-fold.  \vspace{0.15cm} 
       
{\bf Example (2.4)}. As compared with 
Calabi-Yau 3-folds, the statement of (2.2) 
for symplectic varieties is quite simple. 
We shall briefly discuss the difference by 
comparing some examples. 

Let us consider the situation where $\tilde 
X$ is a smooth Calabi-Yau 3-fold containing 
a smooth divisor $E$ and $\pi : 
\tilde X \to X$ is a projective birational 
contraction map of $E$ to a curve $C \subset X$. 
We assume that $C$ is a smooth curve and 
$E \to C$ is a conic bundle with no multiple 
fibers. Denote by $g$ the genus of $C$ 
and denote by $n$ the number of the singular 
fibers of the conic bundle. 
Such birational contractions are 
studied in [Wil, Gr, Na 4].  

$X$ has $A_1$ singularities along $C$; 
hence $\Sigma = {\mathrm{Sing}}(X)$ is 
isomorphic to $C$. Corresponding to 
$n$ singular fibers, $X$ has exactly 
$n$ dissident points.  

First note that any $g$ is possible; hence 
$\omega_C$ is not necessarily trivial ( 
Compare with the symplectic case (1.6)). 
Moreover, the dissident locus $\Sigma_0$ 
has codimension 3. 
Proposition (2.1) is no more true; ${\bold T}^1_X$ 
and ${\bold T}^1_U$ are isomorphic, but 
$H^1(\pi^{-1}(U), \Theta_{\pi^{-1}(U)})$ is 
an infinitesimal dimensional $\bold C$ 
vector space.  
  
Assume that $g = 0$ (i.e. $C = {\bold P}^1$) 
and $n \geq 3$. Then ${\mathrm{Def}}(\tilde X)$ 
and ${\mathrm{Def}}(X)$ are both smooth, but 
$\dim {\mathrm{Def}}(X) = \dim 
{\mathrm{Def}}(\tilde X) + 2n-2-k$, where 
$k := b_2(\tilde X) - b_2(X)$. By the natural 
map $\pi_*$, ${\mathrm{Def}}(\tilde X)$ is 
embedded into ${\mathrm{Def}}(X)$. However, 
$\dim {\mathrm{Def}}(X) > \dim 
{\mathrm{Def}}(\tilde X)$ if $n \geq 4$.    

When $g = 1$ and $n = 0$, we are in a similar 
situation to (2.2), that is, ${\mathrm{Def}}(\tilde X)$ 
and ${\mathrm{Def}}(X)$ are both smooth and 
$\pi_*$ is a finite covering. 
\vspace{0.2cm}

In Theorem (2.2) we have studied a singular 
variety $X$ which has a symplectic resolution. 
However we must often deal with a {\em symplectic} 
variety which does not have a symplectic 
resolution; for example, such varieties appeared 
in [O] as the moduli spaces of rank 2 
semi-stable sheaves on a K3 surface with 
$c_1 = 0$ and with even $c_2 \geq 6$. 
(When $c_1 = 0$ and $c_2 = 4$, the moduli space 
has a 10 dimensional symplectic resolution and 
it provides us with a new example of a 
symplectic manifold.)  Finally we shall 
prove an unobstructedness result for such 
singular symplectic varieties.  \vspace{0.2cm}

{\bf Theorem (2.5)}. {\em Let $X$ be a projective 
symplectic variety. 
Let $\Sigma \subset X$ be the singular locus.  
Assume that $\mathrm{codim}(\Sigma \subset X) \geq 4$. 
Then ${\mathrm{Def}}(X)$ is smooth. }  \vspace{0.15cm}

{\bf Remark}. When $X$ has a symplectic resolution 
$\pi : \tilde X \to X$, the result easily follows 
from Proposition (2.1) because 
${\mathrm{Def}}(\tilde X)$ 
is smooth by a theorem of Bogomolov. \vspace{0.15cm}

{\em Proof}. The following result of Ohsawa 
[Oh] is a key. We shall give an algebraic proof.
\vspace{0.15cm}

{\bf Lemma (2.6)}. {\em The Hodge spectral sequence} 

$$ E_1^{p,q} = H^q(U, \Omega^p_U) =>  
H^{p+q}(U, {\bold C}) $$

{\em degenerates at $E_1$-terms with $p + q = 2$}. 
\vspace{0.15cm}

{\em Proof of Lemma (2.6)}. 
Let $f : Y \to X$ be a resolution of singularities 
such that $E := f^{-1}(\Sigma)$ is a divisor 
with normal crossings and $f^{-1}(U) \cong U$. 
There are natural maps $\phi_{p,q} : 
H^q(Y, \Omega^p_Y(\log E)) 
\to H^q(U, \Omega^p_U)$ induced by the restriction. 
By the mixed Hodge structure on $H^{*}(U, {\bold C})$, 
the spectral sequence 

$$ E^1_{p,q} := H^q(Y, \Omega^p_Y(\log E)) => 
H^{p+q}(Y, \Omega^{\cdot}_Y(\log E)) = 
H^{p+q}(U, {\bold C})$$ 
degenerates at $E_1$ terms. 
Therefore we only have to show that 
$\phi_{p,q}$ are isomorphisms for $p$, $q$ with 
$p + q = 2$. 

By Appendix we will see that $H^q(U, \Omega^{p, alg}_U) 
\cong H^q(U, \Omega^p_U)$ when $p + q = 2$.  
Thus, we have to prove that $\phi^{alg}_{p,q}: 
H^q(Y, \Omega^{p, alg}_Y(\log E)) \to  
H^q(U, \Omega^{p, alg}_U)$ are isomorphisms 
when $p + q = 2$. By the local cohomology sequences, 
it is enough to show that 

(1)  $H^i_E(Y, {\cal O}^{alg}_Y) = 0$ for $i = 2, 3$, 
\vspace{0.12cm}

(2)  $H^i_E(Y, \Omega^{1, alg}_Y(\log E)) = 0$ for 
$i = 1, 2$, and \vspace{0.12cm}

(3)  $H^i_E(Y, \Omega^{2, alg}_Y(\log E)) = 0$ for 
$i = 0, 1$.   \vspace{0.12cm}

For (1) let us show that 
$H^3_E(Y, {\cal O}^{alg}_Y) = 0$.
Let $Y_{\Sigma}$ be the formal completion of 
$Y$ along $f^{-1}(\Sigma)$ and $X_{\Sigma}$ 
the formal completion of $X$ along $\Sigma$. 
$H^3_E(Y, {\cal O}^{alg}_Y)$ is dual to 
$H^{n-3}(Y_{\Sigma}, \omega^{alg}_{Y_{\Sigma}})$. 
Note that 
$H^{n-3}(X_{\Sigma}, 
f_*\omega^{alg}_{Y_{\Sigma}}) = 0$ 
because $\mathrm{codim}(\Sigma \subset X) \geq 4$.
By the Grauert-Riemenschneider vanishing theorem 
and GAGA principle, 
$R^if_*\omega^{alg}_{Y_{\Sigma}} = 0$ for $i >0$. 
Therefore 
$H^{n-3}(Y_{\Sigma}, 
\omega^{alg}_{Y_{\Sigma}}) = 0$. 
The proof that 
$H^2_E(Y, {\cal O}^{alg}_Y) = 0$ is 
similar.   

Next we shall prove (2).  

$H^2_E(Y, \Omega^{1, alg}_Y(\log E))$ is dual to 
$H^{n-2}(Y_{\Sigma}, 
\Omega^{n-1, alg}_Y(\log E)(-E))$. 

Since $\mathrm{Codim}(\Sigma \subset X) 
\geq 4$, the vector spaces $H^{n-2}(X_{\Sigma}, 
f_*\Omega^{n-1, alg}_Y(\log E)(-E))$ and   
$H^{n-3}(X_{\Sigma}, 
R^1f_*\Omega^{n-1, alg}_Y(\log E)(-E)) = 0$ 
are both zero.
 
On the other hand, 
$R^if_*\Omega^{n-1, alg}_Y(\log E)(-E) = 0$ for 
$i \geq 2$ by [St 2] and GAGA principle. Therefore 
$H^{n-2}(Y_{\Sigma}, 
\Omega^{n-1, alg}_Y(\log E)(-E)) 
= 0$.   
We can show that 
$H^1_E(Y, \Omega^{1, alg}_Y(\log E)) = 0$ 
in a similar way. 

The proof of (3) is similar to that of (2).  
\vspace{0.15cm}

{\bf Lemma (2.7)}. {\em Let $X_m$ be an infinitesimal 
deformation of $X$ over $S_m := \mathrm{Spec}(A_m)$, 
where $A_m = {\bold C} [t]/(t^{m+1})$. 
Let $U_m := X_m\vert_U$. 
Then the Hodge spectral sequence} 

$$  E_1^{p,q} = H^q(U, \Omega^p_{U_m/S_m}) 
=>  H^{p+q}(U, A_m)  $$

{\em degenerates at $E_1$-terms with $p + q = 2$. 
In particular, the natural map 
$H^q(U, \Omega^p_{U_m/S_m}) \to 
H^q(U, \Omega^p_{U_{m-1}/S_{m-1}})$ is surjective, 
where $U_{m-1} := U_m \times_{S_m}S_{m-1}$.}  
\vspace{0.2cm}

{\em Proof}. Note that $H^{p+q}(U, A_m) \cong 
H^{p+q}(U, {\bold C})\otimes_{\bold C}A_m $. 
When $m = 0$ the Hodge spectral sequence degenerates 
at $E_1$-terms with $p + q = 2$ by Lemma (2.6). 
Hence 
$\Sigma_{p + q = 2}
\dim_{\bold C}H^q(U, \Omega^p_{U_m/S_m}) = 
\dim_{\bold C}H^2(U, A_m) $. From this it follows 
that $H^q(U, \Omega^p_{U_m/S_m})$ are free $A_m$ 
modules for $p$, $q$ with $p + q = 2$ and 
the Hodge sectral sequence degenerates at 
$E_1$-terms with $p + q = 2$. 
\hspace{0.4cm} Q.E.D.  \vspace{0.2cm}

By Lemma (2.7) the non-degenerate holomorphic 2-form 
$\omega$ extends to an relative 2-form 
$\omega_m \in H^0(U, \Omega^2_{U_m/S_m})$, 
by which $\Theta_{U_m/S_m}$ and 
$\Omega^1_{U_m/S_m}$ are identified. 
Now again by Lemma (2.7) the natural map 
$H^1(U, \Theta_{U_m/S_m}) \to 
H^1(U, \Theta_{U_{m-1}/S_{m-1}})$ is surjective.  
By Proposition (2.1) this implies that the 
$T^1$-lifting property holds for an 
infinitesimal deformation of $X$. 
Therefore ${\mathrm{Def}}(X)$ is smooth.  
\hspace{0.4cm} Q.E.D.   \vspace{0.15cm}

\begin{center}
{\bf 3. Appendix}.       
\end{center}

We shall prove two comparison theorems. 
\vspace{0.12cm}

{\bf Lemma (3.1)}. (1) {\em Let 
$\pi: \tilde X \to X$ be a birational 
projective morphism from a smooth 
symplectic $n$-fold to a normal projective 
variety $X$. Let $U$ be the same as (1.8) 
and put $\tilde U := \pi^{-1}(U)$. 
Then} 

$$ H^1(\tilde U, \Theta^{alg}_{\tilde U}) 
\cong H^1(\tilde U, \Theta_{\tilde U}). $$ 

(2) {\em Let $X$ be a projective 
variety with $\mathrm{Codim}(\Sigma \subset 
X) \geq 4$, where $\Sigma$ is the singular 
locus of $X$. Let $U$ be the regular part 
of $X$. Then} 

$$ H^q(U, \Omega^{p, alg}_U) \cong 
H^q(U, \Omega^p_U) $$  
{\em for $p$ and $q$ with $p + q = 2$.}
\vspace{0.12cm}

{\em Proof}. (1): We shall use 
the same notation as (1.8).
There are two 
exact sequences in the algebraic/ anlytic 
category:

$$0 \to H^1(\Theta^{alg}_U) \to 
H^1(\Theta^{alg}_{\tilde U}) \to
H^0(U, R^1(\pi_U)_*\Theta^{alg}_{\tilde U}) \to 
H^2(\Theta^{alg}_{\tilde U}) $$
and 

$$0 \to H^1(U, \Theta_U) \to H^1(\Theta_{\tilde 
U}) \to H^0(U, R^1(\pi_U)_*\Theta_{\tilde U}) 
\to H^2(\Theta_{\tilde U}). $$

There are natural maps from the first 
sequence to the second one so that they 
make a commutative diagram.  Call these maps 
$\alpha_1$, $\alpha_2$, $\alpha_3$ and $\alpha_4$ 
from the left to the right. We want to prove that 
$\alpha_2$ is an isomorphism. 
  
By a symplectic 2-form $\omega$, one 
has an isomorphism $\Theta^{alg}_{\tilde U} 
\cong \Omega^{1, alg}_{\tilde U}$ (resp. 
$\Theta_{\tilde U} \cong \Omega^1_{\tilde U}$).

By taking the direct image of both sides, 
one has ${\hat\Omega}^1_U \cong 
\Theta_U$, where $\hat{}$ means the 
double dual (cf. [St 1]). 
Since $U$ has only quotient 
singularities, $\mathrm{depth}
({\hat\Omega}^1_U)_p = n$ for 
$p \in U$. Therefore, 
$\mathrm{depth}
(\Theta_U)_p = n$.    
Now, since $\mathrm{Codim}(\Sigma_0 
\subset X) \geq 4$ where $\Sigma_0 
= X - U$, the maps $\alpha_1$ 
and $\alpha_4$ are isomorphisms. 
(cf. [Ha, VI, Theorem 2.1, (a)])\footnote
{The statement of [Ha, VI, Theorem 2.1, (a)] 
is not enough for our purpose. But, as it is 
clear from the proof, the theorem holds under  
certain depth condition.}  

We shall prove that $\alpha_3$ is 
an isomorphism. 
  
At (i) of the proof of Claim 1 in 
Theorem (2.2), we have shown 
that $R^1{\pi_U}*\Omega^1_{\tilde U} 
\cong R^1{\pi_U}\Omega^1_{D_U}$. 
In particular, we see that $R^1{\pi_U}_*
\Omega^1_{\tilde U}$ is a locally 
free ${\cal O}_{\Sigma^{(0)}}$ module, hence 
$R^1{\pi_U}_*\Theta_{\tilde U}$ is 
also a locally free ${\cal O}_{\Sigma^{(0)}}$ 
module. 
$\Sigma^{(0)}$ is a Zariski open set 
of $\Sigma$ and the complement 
$\Sigma - \Sigma^{(0)}$ has 
codmension at least 2 by (1.6). 

Since $R^1\pi_*\Theta^{alg}_{\tilde X}
\otimes_{{\cal O}^{alg}_X}{\cal O}^{alg}
_{\Sigma}\vert_U = R^1{\pi_U}_*
\Theta^{alg}_{\tilde U}$ and  
$R^1\pi_*\Theta_{\tilde X}
\otimes_{{\cal O}_X}{\cal O}
_{\Sigma}\vert_U = R^1{\pi_U}_*
\Theta_{\tilde U}$, we conclude 
that $\alpha_3$ is an isomorphism 
by [Ha, VI, Theorem 2.1].     

Now, by the commutative diagram above, 
$\alpha_2$ is an isomorphism, which completes the 
proof of (1).  \vspace{0.12cm}

(2): This easily follows from 
[Ha, VI, Theorem 2.1, (a)]. 
\vspace{0.2cm}

\begin{center}
Department of Mathematics, 
Graduate school of science, 
Osaka University, Toyonaka, 
Osaka 560, Japan 
\end{center}

\end{document}